\documentclass[11pt,oneside,openany]{article}
\usepackage{tikz}
\usepackage{tikz-cd}

\usepackage{amssymb}
\usepackage{amsthm}
\usepackage{amsmath,mathtools}
\usepackage{amsfonts}
\usepackage{latexsym}
\usepackage{graphics}
\usepackage{fancyhdr}
\usepackage{epstopdf}
\usepackage{setspace}
\usepackage{url}
\usepackage[top=1in, bottom=1in, left=.75in, right=.75in]{geometry}
\usepackage{hyperref}

\usepackage[vcentermath,enableskew]{youngtab}

\newcommand{\downmapsto}{\rotatebox[origin=c]{-90}{\ $\longmapsto$\ }\mkern2mu\hspace{0.2cm}}
\newcommand{\upmapsto}{\hspace{0.2cm}\rotatebox[origin=c]{90}{\ $\longmapsto$\ }\mkern2mu}

\usepackage{enumerate}
\usepackage[shortlabels]{enumitem}
\setenumerate{label={\upshape(\roman*)},
  wide=0pt,topsep=0.2cm,labelsep=0.1em,itemsep=0.15cm,leftmargin=*}

\newtheorem{theorem}{Theorem}[section]

\newtheorem{prop}[theorem]{Proposition}
\newtheorem{lemma}[theorem]{Lemma}
\newtheorem{corollary}[theorem]{Corollary}

\theoremstyle{definition}
\newtheorem{definition}[theorem]{Definition}
\newtheorem{example}[theorem]{Example}

\newtheorem{remark}[theorem]{Remark}

\def\R{{\mathbb R}} 
\def\P{{\mathbb P}} 
\def\Z{{\mathbb Z}} 

\newcommand{\FF}{\mathcal{F}}
\newcommand{\LL}{\mathcal{L}}
\newcommand{\set}[1]{\left\{#1\right\}}        
\DeclareMathOperator{\ddd}{depth}
\newcommand{\dep}[1]{\ddd\left(#1\right)}
\newcommand{\vvec}[2]{\begin{pmatrix}#1\\#2\end{pmatrix}}

\title{On Integer Partitions and Continued Fraction Type Algorithms}
\author{Wael Baalbaki\footnote{Fu Foundation School of Engineering and Applied Sciences, Columbia University, Email: wb2393@columbia.edu} \and Claudio Bonanno\footnote{Dipartimento di Matematica, Universit\`a di Pisa, Largo Bruno Pontecorvo 5, I-56127 Pisa, Italy. Email: claudio.bonanno@unipi.it} \and Alessio Del Vigna\footnote{Dipartimento di Matematica, Universit\`a di Pisa, Largo Bruno
  Pontecorvo 5, I-56127 Pisa, Italy. Email: delvigna@mail.dm.unipi.it} \and Thomas Garrity\thanks{Department of Mathematics and Statistics, Williams College, Williamstown, MA 01267, USA. Email: tgarrity@williams.edu} \and Stefano Isola\footnote{Scuola di Scienze e Tecnologie, Universit\`a di Camerino, via Madonna delle Carceri, I-62032 Camerino, Italy. Email: stefano.isola@unicam.it}}

\date{}
\begin{document}
\maketitle

\begin{abstract}
Our goal is to show that the  additive-slow-Farey version of the Triangle map (a type of multi-dimensional continued fraction algorithm) gives us a method for producing a map from the set  integer partitions of a positive number $n$ into itself. 
We start by showing  that the additive-slow-Farey version of the traditional continued fractions algorithm has a natural interpretation as a method for producing integer partitions of a positive number $n$ into two smaller numbers, with multiplicity.  We provide a complete description of how such integer partitions occur and of the conjugation for the corresponding Young shapes via the dynamics of the classical Farey tree. We use the dynamics of the Farey map to get a new formula for $p(2,n)$, the number of ways for partitioning $n$ into two smaller positive integers, with multiplicity. We then turn to the general case, using the the Triangle map to give  a natural map from general   integer partitions of a positive number $n$  to integer partitions of $n$. This map will still be compatible with conjugation of the corresponding Young shapes.  We will close by the observation that it appears few other multi-dimensional continued fraction algorithms can be used to study partitions.  

\end{abstract}


\section{Introduction}
\label{intro}

This paper is an attempt to put two ideas together, that of integer partitions for positive integers and that of regular and multi-dimensional continued fractions. Of course, the theory of partition numbers is one of the richest areas of a lot of mathematics, especially in combinatorics. Continued fractions and their generalizations to the multi-dimensional case are all division algorithms and important to number theory, to dynamical systems as a rich source of examples and a number of other areas.

We have two audiences in mind: the partition community and the dynamical systems people who are interested in regular and multi-dimensional continued fractions.  Thus there will be a bit more exposition than would be usual in a paper.  This is also why we spend so much time on the somewhat special case of partitions into only two parts.

The main idea behind this paper can initially be seen by the relation between the continued fractions expansion of a rational number $m/n$ and the \emph{Farey tree}, a binary tree containing all rational numbers in $(0,1)$. The construction of the Farey tree recalled in Section \ref{farey} can then be used to generate integer partitions of the denominator $n$. Let us consider for example the number $8/19$, which is studied in details in Examples \ref{ex8/19} and \ref{ex8/19-2}. Our method generates the following partitions of 19:
\begin{eqnarray*}
19 &=& 11 + 8 \\
&=& (8 +3) + 8 = 2\cdot 8 + 3 \\
&=& 2\cdot (5+3) + 3 = 2\cdot 5 + 3\cdot 3\\
&=& 2\cdot (3+2) + 3\cdot 3 = 5\cdot 3 + 2\cdot 2\\
&=& 5\cdot (2+1) +2\cdot 2 = 7\cdot 2+ 5\cdot 1\, .
\end{eqnarray*}
Notice that at each step we are simply applying a slow version of the Euclidean algorithm starting from the couple $(8,19)$. In this slow version of the algorithm we subtract the number 8 from 11 just once, and keep doing the same subtraction until possible. Then we pass to the couple $(3,8)$ and continue. We have thus produced five partitions of 19 into two different parts. We refer to Section \ref{partition} for a background on integer partitions.

The first main result of the paper is that given an integer $n$ we generate all its partitions into two coprime different parts with coprime multiplicities by repeating the previous method for all numbers $m/n \in (0,1)$ with $\gcd(m,n)=1$. We can then extend this result to show how to generate all the partitions of an integer $n$ into two different parts. The relation of the generated partitions with the Farey tree also leads to a formula for the number of these partitions. In Theorem  \ref{two sum}, we give a new formula for the number of ways for partition $n$ into two smaller numbers, with multiplicity.  This formula is quite different than the formula of Kim  \cite{Kim 2012} and has a dynamical interpretation in terms of the Farey tree and the related \emph{Farey map}.

In Section \ref{partition}, we set up our notation for partitions. Section \ref{section-farey-map} deals first with the Farey map, which we refer to as the additive-slow-Farey map to remark its role in the generation of the slow version of the additive Euclidean algorithm. In the same section we also recall the construction of the Farey tree and the interpretation of the Farey map via two-by-two matrices in $SL(2, \Z)$.

In Section \ref{2 case} we introduce our method to find integer partitions of a positive number $n$ into two smaller numbers by using the Farey map and prove the main results of the paper. Our method also produces an extremely natural interpretation of the conjugation among partitions. In the theory of integer partition this conjugation is described in terms of the \emph{Young shape} of a partition and of the \emph{flipping} of the shape. In Theorems \ref{thmconj} and \ref{2-palindorme} we show that this conjugation comes out from the properties of the \emph{binary sequence} of a rational number in the Farey tree and of its reversed sequence.

All of this is about the quite special case of partitioning $n$ into two parts.  It is in Section \ref{npartitions} that we start the discussion of the generalization of our method to partitions into many  parts. This is the part of the paper in which we start using a particular multi-dimensional continued fractions algorithm. In particular we use what we call the \emph{additive-slow-Triangle map}, a  $m$-dimensional version of the Farey map (whose two dimensional version was  introduced and studied in \cite{Bonanno- Del Vigna-Munday} in analogy with the two-dimensional version of the Gauss map defined in \cite{Garrity1}). In Section \ref{npartitions} we will see how the Triangle map can be used to produce various partitions of $n$ into $m$ smaller numbers, with multiplicity.  The geometry of the Triangle map and of its domain is  more  complicated than that for the Farey map. We will determine a  description of which such integer partitions occur and give  description of the conjugation for the corresponding Young shapes via the dynamics of the Triangle map. In Section \ref{$T_D$}, we will further extend the triangle map acting on partitions, giving us a way for capturing all possible partitions into orbits of other partitions.   Now   the Triangle map is only one of many possible multi-dimensional continued fractions algorithms that exist.   In Section \ref{otherMCFs}, we will look at two other well-known multi-dimensional continued fractions algorithms (the M\"{o}nkmeyer map and the Cassaigne map) and show, somewhat surprisingly, that neither can be used to study partition numbers.  Further in that section we discuss that there are only a few multi-dimensional continued fractions algorithms that will create orbits of partition numbers.  We close with questions in the last section. 

\section{Background on Partition Numbers}
\label{partition}

There are many sources for background on partition numbers: the
classical text is Andrews \cite{Andrews} and a good introduction is
Andrews and Eriksson \cite{Andrews-Eriksson}. In this section we
recall what we need in the following.

A \emph{partition} of an integer $n\geq 1$ is a non-increasing
sequence of positive integers
\[
    \lambda = (\lambda_1,\,\lambda_2,\,\ldots,\,\lambda_r)
\]
such that $\lambda_1+\cdots+\lambda_r=n$. In this case we shall write
$(\lambda_1,\,\lambda_2,\,\ldots,\,\lambda_r)\vdash n$. The
\emph{partition number} $p(n)$ is the number of partitions of $n$,
that is the number of ways for adding positive numbers together to get
$n$, with order not mattering. For example $p(6) = 11$ because $6$ can
be obtained in the following $11$ different ways:
\begin{alignat*}{4}
    &6    \hspace{2cm} &&4+1+1   \hspace{2cm}   &&3+1+1+1  \hspace{2cm}    &&2+1+1+1+1\\
    &5+1   &&3+3        &&2+2+2        &&1+1+1+1+1+1\\
    &4+2   &&3+2+1      &&2+2+1+1      &&
\end{alignat*}
If a certain $\lambda_i$ is repeated in the sequence, say $k_i$ times,
we collapse the repeated values and use the compact notation
\[
    (n_1^{k_1},\,\ldots,\,n_m^{k_m}) \vdash n
\]
to denote the partition $n_1\cdot k_1+\cdots +n_m\cdot k_m= n$.  We
shall call $n_1,\,\ldots,\,n_m$ the \emph{parts} and
$k_1,\,\ldots,\,k_m$ the \emph{multiplicities} of the partition. Thus
for $n=6$, we can rewrite the above partitions as
\begin{eqnarray*}
    (6) & \vdash & 6 \\
    (5,1) & \vdash & 6 \\
    (4,2) & \vdash & 6 \\
    (4,1^2) & \vdash & 6 \\
    & \vdots & \\
    (2^2, 1^2) & \vdash & 6 \\
    (2, 1^4) & \vdash & 6 \\
    (1^6)& \vdash & 6
\end{eqnarray*}
As we will be acting on a partition $(n_1^{k_1}, \ldots, n_m^{k_m}) $
via matrices, it will be convenient at some point to use the notation
\[
    (n_1^{k_1},\,\ldots,\,n_m^{k_m}) = (n_1,\,\ldots,\,n_m) \times
    [k_1,\,\ldots,\,k_m],
\]
with round brackets for the parts and square brackets for the
multiplicities.

To a given partition $(\lambda_1,\,\dots,\,\lambda_r)$ we associate
the \emph{Young shape}, a diagram with $k$ rows such that the $i$-th
row contains $\lambda_i$ squares. Equivalently, to a partition of the
form $(n_1^{k_1},\,\ldots,\,n_m^{k_m})$ we associate the shape with
$k_1+\cdots+k_m$ rows such that there are $k_1$ rows with $n_1$
squares on top of $k_2$ rows with $n_2$ squares, and so on.

\begin{example}
    For example, the Young shape for $( 5^3, 3^2, 2^1)\vdash 23$ is
    \[
        \yng(5,5,5,3,3,2)
    \]
\end{example}

\noindent We can always flip any such Young shape, turning the rows
into columns, getting a new Young shape which still represents a
partition of the same integer. We shall refer to this new partition as
to the \emph{conjugate partition} and write
\[
    \lambda \sim_{\mathcal C} \mu
\]
to indicate that the partitions $\lambda$ and $\mu$ are conjugate.

\begin{example}
    Flipping the Young shape of the partition
    $( 5^3, 3^2, 2^1)\vdash 23$ of the previous example gives us the
    Young shape
    \[
        \yng(6,6,5,3,3)
    \]
    which represents the conjugate partition
    $(6^2, 5^1, 3^2) \vdash 23$.
\end{example}

 Explicitly we have
\[
    (n_1^{k_1},n_2^{k_2}) \sim_{\mathcal C}
    ((k_1+k_2)^{n_2},\,k_1^{n_1-n_2}),
\]

\[
    (n_1^{k_1}, n_2^{k_2}, n_3^{k_3})\sim_{\mathcal C}
    ((k_1+k_2+k_3)^{n_3},\, (k_1 + k_2)^{n_2-n_3},\, k_1^{n_1-n_2}).
\]
and in general 
\[
    (n_1^{k_1}, \ldots, n_m^{k_m})\sim_{\mathcal C}
    ((k_1+\ldots + k_m)^{n_m},\, (k_1 +\ldots +   k_{m-1})^{n_2-n_3},\ldots,  k_1^{n_1-n_2}).
\]


\section{Preliminaries on the additive-slow-Farey map} \label{section-farey-map}

In this section we first recap some results from the theory of
continued fractions, we recall the definition of the Farey map and the
construction of the Farey tree, to then show the connection between
the two worlds. (Partially reflecting that different mathematical communities work on these maps, they are sometimes called slow version or additive version of the Gauss map, which is why we sometimes call this the additive-slow-Farey map).

The basics of  continued fractions are in most beginning number theory books.  For more in depth treatment, there is the classic work of Khinchin \cite{khi}.  To see how continued fractions are naturally linked to dynamical systems, see Dajani and Kraakamp \cite{Dajani-Kraakamp} or Hensley \cite{Hensley}.

\subsection{Basic properties of continued fractions} \label{sec:basic-cf}

For $x\in (0,1)$ we denote by $[a_1,\,a_2,\,\ldots]$ its continued
fraction expansion, that is
\[
    x = \frac{1}{a_1+\dfrac{1}{a_2+\dfrac{1}{a_3+\cdots}}},
\]
where $a_j\geq 1$ for all $j$. The expansion is finite if and only if
$x$ is rational and it is also unique provided that for finite
expansions such as $[a_1,\,\ldots,\,a_k]$ we require   $a_k>1$. The
\emph{convergents} of a real number $x\in (0,1)$ are the elements of
the sequence $(p_j/q_j)_{j\geq 0}$ recursively defined as follows:
\begin{alignat*}{3}
    &p_0=0\quad               &\hspace{1.5cm}&q_0=1\\
    &p_1=1\quad               &&q_1=a_1\\
    &p_{j+1} = a_{j+1}p_j+p_{j-1} &&q_{j+1} = a_{j+1}q_j+q_{j-1}.
\end{alignat*}
It is easy to show that
\[
    \frac{p_j}{q_j} = [a_1,\,\ldots,\,a_j].
\]
Note that the sequence of convergents is finite or infinite according
to whether $x$ is rational or not. In particular, if
$x=[a_1,\,\ldots,\,a_k]$ is rational then the sequence of convergents
stops at $\frac{p_k}{q_k}=x$. On the other hand, if $x$ is irrational
then $\frac{p_j}{q_j}\to x$ as $j\to \infty$ and the convergents are
the best rational approximations of $x$ in a precise sense. Another
classical result which we shall use in the following is the so-called
\emph{mirror formula}, that is
\begin{equation}\label{mirror}
    \frac{q_{j-1}}{q_j} = [a_j,\,\ldots,\,a_1].
\end{equation}

\subsection{The Farey map and the Farey tree}
\label{farey}

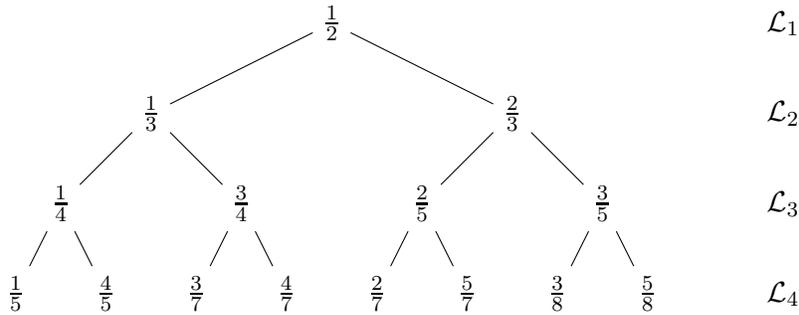
\begin{figure}[h]
    \begin{center}

\begin{tikzpicture}[
    level 1/.style = {sibling distance=4cm},
    level 2/.style = {sibling distance=2cm},
    level 3/.style = {sibling distance=1cm},
    level distance          = 1cm,
    edge from parent/.style = {draw},
    scale=1.2
    ]

    \foreach \n in {1,2,3,4} {
      \pgfmathsetmacro\p{\n}
      \node at (5,-\n+1) {$\mathcal{L}_{\pgfmathprintnumber\p}$};
    }

    \node {$\frac 12$}
    child{
      node {$\frac 13$}
      child{
        node {$\frac 14$}
        child{
          node {$\frac 15$}
        }
        child{
          node {$\frac 45$}
        }
      }
      child{
        node {$\frac 34$}
        child{
          node {$\frac 37$}
        }
        child{
          node {$\frac 47$}
        }
      }
    }
    child{
      node {$\frac 23$}
      child{
        node {$\frac 25$}
        child{
          node {$\frac 27$}
        }
        child{
          node {$\frac 57$}
        }
      }
      child{
        node {$\frac 35$}
        child{
          node {$\frac 38$}
        }
        child{
          node {$\frac 58$}
        }
      }
    };
\end{tikzpicture}

        \caption{The first four levels of the Farey tree.}\label{fig:Farey}
    \end{center}
\end{figure}

Split the unit interval $I=[0,1]$ into  two sub-intervals
$I_0 = \left[\frac{1}{2}, 1\right]$ and
$I_1 = \left[0, \frac{1}{2}\right]$. The \emph{Farey map} is the map
$F:I\to I$ defined to be
\[
    F(x) = \begin{cases}
        F_0(x)=\frac{1-x}{x}\, , &\, \text{if }x \in I_0\\
        F_1(x)=\frac{x}{1-x}\, , &\, \text{if }x \in I_1
    \end{cases}
\]
It has the two local inverses $\Phi_0\coloneqq F_0^{-1}:I\to I_0$
and $\Phi_1\coloneqq F_1^{-1}:I\to I_1$ given by
\[
    \Phi_0(x) = \frac{1}{1+x}\quad\text{and}\quad
    \Phi_1(x)=\frac{x}{1+x}.
\]
It is well-known that the set
$\FF = \bigcup_{k= 0}^\infty F^{-k}(\frac 12)$ contains every rational
number in $(0,1)$ exactly once. Furthermore, $\FF$ can be ordered as a
binary tree, the \emph{Farey tree}. The recursive construction works as
follows: the root of the tree is $\frac 12$ and the two children of
the vertex $\frac pq$ are its backward images under the Farey map,
namely $\Phi_1\big(\frac pq\big)=\frac p{p+q}$ and
$\Phi_0\big(\frac pq\big)=\frac q{p+q}$. Note that each rational
number appears reduced in lowest terms in the tree because the root
$\frac 12$ is. The \emph{levels} of the Farey tree are the sets
$\LL_k=F^{-k+1}(\frac 12)$ for $k\geq 1$, so that
\[
    \LL_1 = \set{\frac 12},\quad
    \LL_2 = \set{\frac 13, \frac 23},\quad
    \LL_3 = \set{\frac 14,\,\frac 25,\,\frac 35,\,\frac
      34},\quad\ldots
\]
See Figure~\ref{fig:Farey} to see the structure of the first few
levels of the Farey tree.

\begin{remark}
    We are describing here the tree constructed with the inverse
    branches of the Farey map, which we call the  Farey tree, with a slight
    abuse. Indeed the ``actual'' Farey tree $\tilde{\FF}$ is the one
    defined by the mediant between neighboring fractions (see
    \cite{BI} for instance). It can be shown that the levels of
    $\tilde{\FF}$ and those of ``our'' Farey tree $\FF$ coincide, but
    in $\tilde{\FF}$ the fractions of each level appear in ascending
    order.
\end{remark}

\noindent Sometimes it is convenient to extend the Farey tree above
the root adding $\frac 11$, which is mapped to $\frac 12$ by both
$\Phi_0$ and $\Phi_1$. In this case we also set
$\LL_0=\set{\frac 11}$.

\begin{definition}
    For $\frac pq\in (0,1)$ in lowest terms we shall call the
    \emph{depth} of $\frac pq$ the level of the Farey tree $\frac pq$
    belongs to. That is, we shall write
    \[
        \dep{\frac pq}=k \quad\Leftrightarrow\quad \frac pq\in \LL_k.
    \]
    If $\frac pq$ is not in lowest terms, we define its depth as the
    depth of the reduced form of $\frac pq$.
\end{definition}

\noindent For instance $\dep{\frac 23}=\dep{\frac{10}{15}}=2$ and
$\dep{\frac 3{8}}=4$.

Starting from a given $x\in (0,1)$, we can iterate the map $F$, with
the iterations terminating whenever the image is $\frac 12$. In
particular, we can encode a given $x\in (0,1)$ as a sequence of zeros
and ones if we keep track of whether the iterations fall in
$J_0=(\frac{1}{2},1)$ or in $J_1=(0,\frac{1}{2})$. More precisely, we
associate to $x$ the unique binary word
$\sigma (x)=\sigma_1\sigma_2\cdots$ such that
\[
    \sigma_j\in \{0,1\}\quad\text{and}\quad F^{j-1}(x)\in
    J_{\sigma_j}\quad\text{for every $j\geq 1$}.
\]
Of course $\sigma(\frac{1}{2})$ is the empty word. We shall refer to
$\sigma(x)$ as the \emph{binary sequence} of $x$. We know that
$\sigma(x)$ is finite if and only if $x$ is rational: this follows
immediately because the Farey tree is constructed by taking backward
images of $\frac 12$ under $F$ and it contains all and only the
rational numbers in $(0,1)$. Let $x=\frac pq\in (0,1)$ be a rational
number and let $\sigma(x)=\sigma_1\cdots \sigma_\ell$. We have
$\ell=\dep{x}-1$ and, from the definition of $\sigma$, also
\[
    F_{\sigma_\ell}\circ\cdots\circ F_{\sigma_1}\left(\frac pq\right)
    = \frac 12,
\]
which is the same as
\begin{equation}\label{eq:rational}
    \frac pq = \Phi_{\sigma_1}\circ \cdots
    \circ\Phi_{\sigma_\ell}\left(\frac 12\right).
\end{equation}

\begin{example}\label{ex8/19}
    The orbit of $\frac 8{19}$ under $F$ is
    \[
        \frac{8}{19} \overset{F_1}{\longmapsto} \frac{8}{11}
        \overset{F_0}{\longmapsto} \frac{3}{8} \overset{F_1}{\longmapsto}
        \frac{3}{5} \overset{F_0}{\longmapsto} \frac{2}{3}
        \overset{F_0}{\longmapsto} \frac{1}{2},
    \]
    so that $\dep{\frac 8{19}}=6$,
    $\sigma\left(\frac 8{19}\right)=10100$, and
    $\frac 8{19} = \Phi_1\circ\Phi_0\circ\Phi_1\circ\Phi_0\circ\Phi_0
    \left(\frac 12\right)$.
\end{example}

There is a beautiful connection between the Farey tree and continued
fractions. Indeed the Farey map has a simple action on continued
fraction expansions. For $x=[a_1,\,a_2,\,a_3,\,\ldots]$ (the expansion
may be finite or not) we have
\[
    F([a_1,\,a_2,\,a_3,\,\ldots]) =
    \begin{cases}
        [a_1-1,\,a_2,\,a_3,\,\ldots]\, , & \text{if }a_1>1\\
        [a_2,\,a_3,\,\ldots]\, , & \text{if }a_1=1
    \end{cases}.
\]
In other words, $F$ subtracts off a $1$ from the first digit of the
expansion if it is greater than 1 and deletes it when it is $1$. Note that
$a_1>1$ if and only if $x=[a_1,\,a_2,\,a_3,\,\ldots]$ is in the
interval $(0,\frac 12)$, so that more precisely $F_1$ acts subtracting
off a $1$ from $a_1$, while $F_0$ acts by deleting it. We now consider
again the case $x=\frac pq$ rational to see
Equation~\eqref{eq:rational} in a new light. We know that
$x=[a_1,\,\ldots,\,a_k]$ has a finite continued fraction expansion
(which is unique, as long as $a_k>1$) and thus, considering the action
of $F$ on the expansion of $x$, we have
\[
    \sigma(x)= 1^{a_1-1}01^{a_2-1}0\cdots 1^{a_k-2}.
\]
We remark that if $k=1$ we have $\sigma(x)=1^{a_1-2}$. Hence
Equation~\eqref{eq:rational} can be rewritten as
\begin{align}\label{eq:rational2}
    \frac pq &= \Phi_1^{a_1-1}\Phi_0\circ
    \cdots\circ\Phi_1^{a_{n-1}-1}\Phi_0 \circ
    \Phi_1^{a_k-2}\left(\frac 12\right)=
    \nonumber\\
    &= \Phi_1^{a_1-1}\Phi_0\circ
    \cdots\circ\Phi_1^{a_{k-1}-1}\Phi_0 \circ
    \Phi_1^{a_k-1}\left(\frac 11\right).
\end{align}
With the same argument we can also write analogous equations for the
convergents $\frac{p_j}{q_j}$. A consequence of
Equation~\eqref{eq:rational2} is
\[
    \dep{\frac pq} = \sum_{j=1}^k a_j -1.
\]
For more details the reader can refer to \cite{BI,I}.

\subsection{Rewriting the Farey map as matrix multiplication}
\label{sec:Fmatrix}

We can also interpret the Farey map as not acting on the unit interval
$I$ in $\R$ but as acting on the cone
\[
    \set{ \vvec{x}{y} \in \R^2 \,:\, y \geq x \geq
      0}\setminus\set{\vvec{0}{0}}.
\]
Now the vector $(x,y)$ is representing the real number $\frac xy$, so
that the unit interval is actually in bijection with the lines in the
cone. In particular, we represent the rational number $\frac pq$ with
the vector $\begin{psmallmatrix}p\\q\end{psmallmatrix}$. Hence the
action of $F$ becomes
\[
    F\vvec{x}{y} =
    \left\{\begin{array}{cc}
        F_0\vvec{x}{y} = \vvec{y-x}{x}, & \text{if $x\geq y-x$}\\[0.5cm]
        F_1\vvec{x}{y}= \vvec{x}{y-x}, & \text{if $y-x\geq x$}
    \end{array}\right..
\]
In this way the Farey map can be seen as the action by left
multiplication of column vectors by $2\times 2$ matrices, that is
\[
    F_0\vvec{x}{y} = \begin{pmatrix} -1&1\\1&0\end{pmatrix}
    \vvec{x}{y} = \vvec{y-x}{x}
    \quad\text{and}\quad
    F_1\vvec{x}{y} = \begin{pmatrix} 1&0\\-1&1\end{pmatrix}
        \vvec{x}{y} = \vvec{x}{y-x}.
\]
We highlight the two matrices
\[
    F_0=\begin{pmatrix} -1&1\\1&0\end{pmatrix}
    \quad\text{and}\quad
    F_1= \begin{pmatrix} 1&0\\-1&1\end{pmatrix}.
\]
The two inverse branches of $F$ can also be expressed in terms of
multiplication of two matrices, which we denote by
\[
    \Phi_0 = \begin{pmatrix} 0 & 1\\ 1 & 1 \end{pmatrix}
    \quad\text{and}\quad
    \Phi_1 = \begin{pmatrix} 1 & 0 \\ 1 & 1 \end{pmatrix}.
\]
With this convention, Equation~\eqref{eq:rational} can be rewritten as
\[
    \vvec{p}{q} = \prod_{j=1}^{\ell} \Phi_{\sigma_j} \vvec{1}{2}
    = \prod_{j=1}^{\ell+1} \Phi_{\sigma_j} \vvec{1}{1},
\]
where $\ell=\dep{p/q}-1$ and the last digit $\sigma_{\ell+1}$ can be
either $0$ or $1$. Thus each rational $\frac pq$ is encoded by two
matrices, namely
\[
    \prod_{j=1}^{\ell} \Phi_{\sigma_j}\cdot \Phi_0
    \quad\text{and}\quad
    \prod_{j=1}^{\ell} \Phi_{\sigma_j}\cdot \Phi_1.
\]
Note that they have the same columns in the two possible orders, so
that in particular one of the two has determinant $1$ and the other
one has determinant $-1$. We shall refer to the matrix with positive
determinant as \emph{the} matrix of $\frac pq$. We now highlight an
important property of these matrices, which we shall use in the proof
of Theorem~\ref{2parts-farey}.

\begin{lemma}
    \label{lemma:matrixpq}
    The set of the matrices of the rationals in the interval $(0,1)$
    is
    \[
        \set{\begin{pmatrix} p'&p''\\q'&q''\end{pmatrix}\in M(2,\Z)
          \,:\,1\leq p'\leq q',\ 0\leq p''<q'',\ p'q''-p''q'=1}\subseteq
        SL(2,\Z).
    \]
\end{lemma}

\begin{proof}
    Let $\frac pq$ be a rational in the interval $(0,1)$ and let
    \[
        \begin{pmatrix} p'&p''\\q'&q''\end{pmatrix}
    \]
    be its matrix. By definition $p'q''-p''q'=1$ and an easy induction
    shows that $1\leq p'\leq q'$ and $0\leq p''<q''$. One may also
    look at the two columns as representing the two rationals
    $\frac{p'}{q'}<\frac{p''}{q''}$ with
    $\frac pq=\frac{p'+p''}{q'+q''}$. As it is well-known the
    condition $p'q''-p''q'=1$ means that the two rationals are
    neighbours in the Farey sequence and that the matrix of $\frac pq$
    is nothing but the matrix of the coding described in \cite{BI},
    and the thesis now follows.
\end{proof}

For instance, Example~\ref{ex8/19} can be rewritten in terms of
matrices as
\[
    \vvec{8}{19} = \Phi_1\Phi_0\Phi_1\Phi_0\Phi_0 \vvec{1}{2},
\]
and the matrix of $\frac 8{19}$ is
\[
    \Phi_1\Phi_0\Phi_1\Phi_0\Phi_0\Phi_0=\begin{pmatrix}
        3&5\\7&12\end{pmatrix}.
\]


\section{Partitions into two different parts}\label{2 case}

\subsection{Partitions generated by the Farey map}
\label{2-intro}

In this section we show how to generate integer partitions into two
parts by using the dynamics of the Farey map and the Farey tree.

Let $n\geq 2$ be an integer and $r$ be such that $1\leq r< n$ and
$(r,n)=1$, so that $\frac rn$ appears in some level of the Farey
tree. If $\sigma(\frac rn)=\sigma_1\cdots\sigma_\ell$ is the binary
sequence of $\frac rn$ introduced in Section~\ref{farey}, we have
\[
    \vvec{r}{n} = \prod_{j=1}^\ell \Phi_{\sigma_j} \vvec{1}{2} =
    \prod_{j=1}^{\ell+1} \Phi_{\sigma_j} \vvec{1}{1},
\]
where $\sigma_{\ell+1}$ can be either $0$ or $1$. Recall that
$\ell=\dep{r/n}-1$. For each $m=0,\,\ldots,\,\ell+1$ we split the
above product as
\begin{equation}
    \label{generate-part}
    \vvec{r}{n} = \prod_{j=1}^m \Phi_{\sigma_j} \cdot
    \prod_{j=m+1}^{\ell+1} \Phi_{\sigma_j} \vvec{1}{1}
\end{equation}
and we set
\[
    \begin{pmatrix} h_2(m)&h_1(m)\\k_2(m)&k_1(m)\end{pmatrix} =
    \prod_{j=1}^m \Phi_{\sigma_j}
    \quad\text{and}\quad \vvec{n_2(m)}{n_1(m)} =
    \prod_{j=m+1}^{\ell+1} \Phi_{\sigma_j} \vvec{1}{1}.
\]
Explicitly writing the second component of Equation~\eqref{generate-part}, we get
\[
    n = k_1(m)n_1(m) + k_2(m)n_2(m),
\]
with $n_1(m)\geq n_2(m)$. Thus for each $m$ we get a partition of $n$,
namely $(n_1(m)^{k_1(m)}, n_2(m)^{k_2(m)})\vdash n$, and the integer
$m$ is called the \emph{generation} of the partition. Note that for
$m=0$ Equation~\eqref{generate-part} reads
\[
    \vvec{r}{n} = \begin{pmatrix}1&0\\0&1\end{pmatrix} \vvec{r}{n},
\]
which induces the partition $(n^1,r^0)=(n^1)$ of generation $m=0$. On
the other hand for $m=\ell+1$, whether $\sigma_{\ell+1}=0$ or
$\sigma_{\ell+1}=1$, we have $k_1(m)+k_2(m)=n$ and thus the induced
partition is $(1^{k_2(m)},1^{k_1(m)})=(1^n)$. Since we are interested
in partitions into two distinct parts here, we do not include these
two cases in the following definition.

\begin{definition}
    We shall call the sequence of partitions
    $(n_1(m)^{k_1(m)}, n_2(m)^{k_2(m)})$ for $m=1,\,\ldots,\,\ell$ the
    \emph{orbit of partitions generated by $\frac rn$}.
\end{definition}

\begin{example} \label{ex8/19-2}
    We again consider Example~\ref{ex8/19} to see how the above
    construction works in practice. Set $n=19$ and $r=8$, so that
    \[
        \vvec{8}{19} = \Phi_1\Phi_0\Phi_1\Phi_0\Phi_0\vvec{1}{2}=
        \Phi_1\Phi_0\Phi_1\Phi_0\Phi_0\Phi_{i}\vvec{1}{1}
    \]
    with $i=0$ or $i=1$. Choosing $i=0$ we then get the following
    partitions:
    \begin{alignat*}{3}
        m&=0: \hspace{1cm} &&\vvec{8}{19} = I\vvec{8}{19}
        \hspace{6.5cm}&&(8^0,19^1)=(19^1)\\
        m&=1:  &&\vvec{8}{19} = \Phi_1\vvec{8}{11} = \begin{pmatrix} 1 & 0 \\ 1 & 1 \end{pmatrix}\vvec{8}{11}
        &&(11^1,8^1)\\
        m&=2:  &&\vvec{8}{19} = \Phi_1\Phi_0\vvec{3}{8} = \begin{pmatrix} 0 & 1 \\ 1 & 2 \end{pmatrix}\vvec{3}{8}
        &&(8^2,3^1)\\
        m&=3:  &&\vvec{8}{19} = \Phi_1\Phi_0\Phi_1\vvec{3}{5} = \begin{pmatrix} 1 & 1 \\ 3 & 2 \end{pmatrix}\vvec{3}{5}
        &&(5^2,3^3)\\
        m&=4:  &&\vvec{8}{19} = \Phi_1\Phi_0\Phi_1\Phi_0\vvec{2}{3} = \begin{pmatrix} 1 & 2 \\ 2 & 5 \end{pmatrix}\vvec{2}{3}
        &&(3^5,2^2)\\
        m&=5:  &&\vvec{8}{19} = \Phi_1\Phi_0\Phi_1\Phi_0\Phi_0\vvec{1}{2} = \begin{pmatrix} 2 & 3 \\ 5 & 7 \end{pmatrix}\vvec{1}{2}
        &&(2^7,1^5)\\
        m&=6:  &&\vvec{8}{19} = \Phi_1\Phi_0\Phi_1\Phi_0\Phi_0\Phi_0\vvec{1}{1} = \begin{pmatrix} 3 & 5 \\ 7 & 12 \end{pmatrix}\vvec{1}{1}
        &&(1^{12},1^7)=(1^{19})
    \end{alignat*}
    In case $i=1$ the last partition which we get is
    $(1^7,1^{12})=(1^{19})$.
\end{example}

Note that the orbit of partitions of a fraction contains pairwise
distinct partitions. But different fractions can induce orbits of
partitions sharing some element and can even induce the same orbit of
partitions, as the following lemma shows.

\begin{lemma}
    \label{lemma:paired}
    Let $n\geq 2$ be an integer and $r$ be such that $1\leq r< n$ and
    $(r,n)=1$. Then $\frac rn$ and $\frac{n-r}n$ induce the same orbit
    of partitions of $n$.
\end{lemma}
\begin{proof}
    For $n=2$ the statement is trivial, so we let $n>2$. Note that the
    fractions $\frac rn$ and $\frac{n-r}n$ are the two children of
    $\frac r{n-r}$ in the Farey tree, thus they have the same depth
    and their forward orbits under $F$ coincide. Without loss of
    generality we assume $r<\frac n2$, so that
    $\frac rn<\frac 12<\frac{n-r}r$ and
    \[
        \vvec{r}{n} = \Phi_1\vvec{r}{n-r}
        \quad\text{and}\quad
        \vvec{n-r}{n} = \Phi_0\vvec{r}{n-r}.
    \]
    Moreover, their binary sequences are
    $\sigma(\frac rn)=1\sigma_2\cdots\sigma_\ell$ and
    $\sigma(\frac{n-r}n)= 0\sigma_2\cdots\sigma_\ell$. As a
    consequence, the induced partitions at generation $m=1$ are both
    equal to $((n-r)^1,r^1)$. Now if $\ell\geq 2$, let
    $m=2,\,\ldots,\,\ell$ and consider
    \[
        \vvec{r}{n} = \Phi_1\prod_{j=2}^m\Phi_{\sigma_j} \vvec{n_2(m)}{n_1(m)}
        \quad\text{and}\quad
        \vvec{n-r}{n} = \Phi_0\prod_{j=2}^m\Phi_{\sigma_j}\vvec{n_2(m)}{n_1(m)}.
    \]
    It is now straightforward to prove that the bottom rows of the two
    matrices $\Phi_1\prod_{j=2}^m\Phi_{\sigma_j}$ and
    $\Phi_0\prod_{j=2}^m\Phi_{\sigma_j}$ coincide, and so do the
    induced partitions of $n$.
\end{proof}

We now characterise precisely the partitions which can be obtained
through the Farey tree and the Farey map. We shall give a second proof
of this result in Section~\ref{p(2,n)}.

\begin{theorem}
    \label{2parts-farey}
    Let $n\geq 2$ be an integer. A partition
    $(n_1^{k_1},n_2^{k_2})\vdash n$ can be obtained from the dynamics
    of the Farey map if and only if $(n_1,n_2)=1$ and
    $(k_1,k_2)=1$.
\end{theorem}

\begin{proof}
    $(\Rightarrow)$ If $(n_1,n_2)>1$ then the fraction
    $\frac{n_2}{n_1}$ is not in lowest terms and thus it does not
    appear on the Farey tree. If $(k_1,k_2)>1$ then a matrix of the form
    \[
        \begin{pmatrix} *&*\\k_2&k_1\end{pmatrix}
    \]
    cannot have determinant $\pm 1$, thus it cannot be a finite
    product of the matrices $\Phi_0$ and $\Phi_1$.
    \\[0.15cm]
    $(\Leftarrow)$ We are now given a partition
    $(n_1^{k_1},n_2^{k_2})\vdash n$ with $(n_1,n_2)=1$ and
    $(k_1,k_2)=1$. To prove that it is induced from the Farey tree it
    suffices to show that there exist two integers $h_1$ and $h_2$
    such that $0\leq h_1< k_1$, $1\leq h_2\leq k_2$, and
    $h_2k_1-h_1k_2=1$. Indeed, if this holds then
    Lemma~\ref{lemma:matrixpq} shows that the matrix
    \[
        \begin{pmatrix} h_2&h_1\\k_2&k_1\end{pmatrix}
    \]
    is a finite product of the matrices $\Phi_0$ and $\Phi_1$, thus
    by setting
    \[
        \vvec{r}{n} =
        \begin{pmatrix} h_2&h_1\\k_2&k_1\end{pmatrix}\vvec{n_2}{n_1}
    \]
    we have that $\frac rn$ is a fraction on the Farey tree and
    the partition $(n_1^{k_1},n_2^{k_2})$ is a member of the orbit of partitions generated by $\frac rn$.
    \\[0.15cm]
    Since $(k_1,k_2)=1$ there exist two integers $\tilde h_1$ and
    $\tilde h_2$ such that $\tilde h_2k_1-\tilde h_1k_2=1$. All the
    solutions to the equation $h_2k_1-h_1k_2=1$ are
    \[
        h_1(t) = \tilde h_1 - t k_1\quad\text{and}\quad
        h_2(t) = \tilde h_2 - t k_2,
    \]
    where $t$ could be any integer number. Without loss of generality
    we can then assume that both $\tilde h_1$ and $\tilde h_2$ are
    positive. We choose $t_*=\left[\frac{\tilde h_1}{k_1}\right]$, so
    that $0\leq h_1(t_*)< k_1$. Then
    \begin{align*}
        h_2(t_*) &= \tilde h_2 - \bigg[\frac{\tilde h_1}{k_1}\bigg] k_2 =
        \tilde h_2 - \left(\frac{\tilde h_1}{k_1}-\bigg\{\frac{\tilde
              h_1}{k_1}\bigg\}\right)k_2=
        \\
        &= \frac{\tilde h_2k_1 - \tilde
          h_1k_2}{k_1}+\bigg\{\frac{\tilde h_1}{k_1}\bigg\}k_2= \frac
        1{k_1}+\bigg\{\frac{\tilde h_1}{k_1}\bigg\}k_2,
    \end{align*}
    so that $1\leq h_2(t_*)\leq k_2$. By setting $h_1=h_1(t_*)$ and
    $h_2=h_2(t_*)$.
\end{proof}

Let $p_F(2,n)$ denote the number of partitions of $n$ into two
different parts obtained by the Farey map. Thanks to the previous
Theorem~\ref{2parts-farey} we can give a formula for $p_F(2,n)$.

\begin{corollary}
    \label{pF(2,n)}
    For $n\geq 2$ we have
    \begin{equation}
        \label{eq:pF(2,n)}
        p_F(2,n) = \frac 12 \left(\sum_{\substack{r=1\\(r,n)=1}}^{n-1}
        \dep{\frac rn} - \varphi(n)\right),
    \end{equation}
    where $\varphi(n)$ is the Euler totient function. Moreover,
    $p_F(2,n)=p(2,n)$ if and only if $n$ is prime or $n=4$.
\end{corollary}

\begin{proof}
    Every $\frac rn$ in the Farey tree, that is with $1\leq r<n$ and
    $(r,n)=1$, generates an orbit of partitions of $n$. Each of these
    orbits contains $\dep{\frac rn}-1$ pairwise distinct partitions
    because there are just as many fractions strictly above $\frac rn$
    in the Farey tree and up to $\frac 12$. Counting of all of
    these partitions yields
    \[
        \sum_{\substack{r=1\\(r,n)=1}}^{n-1} \left(\dep{\frac rn}
          -1\right) = \sum_{\substack{r=1\\(r,n)=1}}^{n-1} \dep{\frac
          rn} - \varphi(n).
    \]
    Lemma~\ref{lemma:paired} shows that each $\frac rn$ is paired with
    $\frac{n-r}{n}$, in the sense the they induce the same orbit of
    partitions. Moreover, the orbit of partitions generated by two
    non-paired fractions are disjoint. Thus in the above counting each
    partition appears exactly twice.
    \\[0.15cm]
    As for the second part of the theorem we prove the two
    implications separately.
    \\[0.15cm]
    $(\Rightarrow)$ If $n$ is prime and
    $(n_1^{k_1},n_2^{k_2})\vdash n$ then necessarily $(n_1,n_2)=1$ and
    $(k_1,k_2)=1$, and thus by Theorem~\ref{2parts-farey} the
    partition $(n_1^{k_1},n_2^{k_2})$ can be obtained from the
    dynamics of the Farey map. The case $n=4$ can be verified
    explicitly.
    \\[0.15cm]
    $(\Leftarrow)$ Suppose that $n$ is composite and $n\neq 4$, so
    that $n=ab$ for some $a>b\geq 2$. If $b\geq 3$ then the partition
    $((b-1)^a,1^a)$ cannot be obtained from the dynamics of the Farey
    map because the multiplicities are not relatively prime. If $b=2$
    then $n=2a$, $a\geq 3$, and the partition $((2(a-1))^1,2^1)$
    cannot be obtained as well.
\end{proof}

\begin{remark}
    Corollary~\ref{pF(2,n)} can be compared with the following purely
    number theoretical expression obtained by Kim in \cite{Kim
      2012}:
    \begin{equation}
        \label{kim}
        p(2,n)= \frac 12 \left(\sum_{r=1}^{n-1}
          \sigma_0(r)\sigma_0(n-r) - \sigma_1(n)+\sigma_0(n)\right),
    \end{equation}
    where $\sigma_j(n)=\sum_{d\mid n} d^j$. Note that when $n$ is
    prime then $\varphi(n)=n-1$ and $\sigma_1(n)-\sigma_0(n)=n-1$ do
    coincide\footnote{More precisely, when $n>1$ then
      $\sigma_1(n)-\sigma_0(n)=\sum_{d\mid n}(d-1)=n-1 + \sum_{d\mid
        n,\ d<n} (d-1) \geq \varphi(n)$ and the equality holds if and
      only if $n$ is prime.}. However the two sums of
    \eqref{eq:pF(2,n)} and \eqref{kim} does not in general hold
    termwise, \emph{e.g.} take $n=11$ and $r=3$, so that
    $\dep{\frac 3{11}}=5$ but
    $\sigma_0(3)\cdot \sigma_0(8)=2\cdot 4=8$.
\end{remark}

\subsection{Conjugate partitions and continued fractions: Palidromes}
\label{2-conjugation}

Let $n\geq 2$ and consider a partition
$(n_1^{k_1},n_2^{k_2})\vdash n$. As recalled in
Section~\ref{partition}, by flipping its Young shape one obtains a
new partition of $n$, the conjugate partition
$({\tilde n}_1^{{\tilde k}_1}, {\tilde n}_2^{{\tilde k}_2})$, with
${\tilde n}_1=k_1+k_2$, ${\tilde n}_2=k_1$, ${\tilde k}_1=n_2$, and
${\tilde k}_2=n_1-n_2$. We showed that partitions into two parts can
be generated by the dynamics of the Farey map, thus one may wonder
whether there is a way to dynamically characterise also the
conjugacy: this is the content of this section.

Consider again Equation~\eqref{eq:rational}, which express every
rational in the Farey as a backward image of $\frac 12$. Since
$\frac 12 = \Phi_0(\frac 11) = \Phi_1(\frac 11)$ we can also rewrite
the equation as
$\frac rn = \Phi_{\sigma_1}\circ \cdots \circ\Phi_{\sigma_\ell}\circ
\Phi_{\sigma_{\ell+1}}\left(\frac 11\right)$, where $\sigma_{\ell+1}$
can be either $0$ or $1$. Here we make the choice $\sigma_{\ell+1}=1$
and we shall call $\sigma_1\cdots\sigma_{\ell}1$ the \emph{extended
  binary sequence} of $\frac rn$. Note that if $1\leq r <\frac{n}2$
then $\sigma_1=1$.

\begin{theorem}[Palidrome Version 1] \label{thmconj} Let $r$ be such that
    $1\leq r <\frac{n}2$ and $(r,n)=1$, and suppose that $\frac rn$
    has extended binary sequence
    $\sigma_1\cdots \sigma_{\ell}1 = 1\sigma_2\cdots \sigma_{\ell}1$,
    with $\ell=\dep{\frac rn}-1$. Let $(n_1^{k_1},n_2^{k_2})$ be a
    partition of $n$ in the $m$-th generation of the orbit of
    partitions generated by $\frac rn$, that is suppose there is
    $1\leq m \leq \ell$ such that
    \[
        \vvec{r}{n}= \prod_{j=1}^{m}\Phi_{\sigma_j}\vvec{n_2}{n_1}
        = \begin{pmatrix} t & s \\ k_2 & k_1 \end{pmatrix}
        \vvec{n_2}{n_1}.
    \]
    Then the conjugate partition
    $({\tilde n}_1^{{\tilde k}_1}, {\tilde n}_2^{{\tilde k}_2})$ is in
    the $(\ell+1-m)$-th generation of the orbit generated by
    $\frac{\tilde r}{n}$, where $\frac{\tilde r}{n}$ is the fraction
    with extended binary sequence
    $1\sigma_\ell\cdots\sigma_1 = 1\sigma_\ell\cdots\sigma_21$. In
    other words,
    \[
        \vvec{\tilde r}{n}= \prod_{j=1}^{\ell+1-m}
        \Phi_{\sigma_j}\vvec{\tilde n_2}{\tilde n_1} = \begin{pmatrix}
            \tilde t & \tilde s \\ \tilde k_2 & \tilde
            k_1 \end{pmatrix} \vvec{\tilde n_2}{\tilde n_1}.
    \]
    Moreover, the continued fraction expansion of $\frac{\tilde r}n$
    has the same digits of that of $\frac rn$, but in reversed order.
\end{theorem}
\begin{proof}
    As already noted, the extended binary sequence
    $\sigma_1\cdots \sigma_{\ell}1$ of $\frac rn$ starts with
    $\sigma_1=1$. Then by reversing it we have another sequence
    starting by $1$ and it represents a fraction $\frac{\tilde r}{n}$
    in the same level of the Farey tree as $\frac rn$. By
    Equation~\eqref{eq:rational2}, the overall effect of this amounts
    simply to reverse the order of the partial quotients in the
    continued fractions expansion, that is
    \[
        \frac rn =[a_1,\,\ldots,\,a_k] \quad\text{if and only if}\quad
        \frac {\tilde r}n =[a_k,\,\ldots,\,a_1].
    \]
    Note that, since $1\leq r < \frac{n}2$ then $a_1>1$, so that the
    expansion for $\frac {\tilde r}n$ is unambiguously defined, and
    viceversa.  Setting $p_{-1}=q_0 =1$, $q_{-1}=p_0=0$, and
    $\frac {p_j}{q_j} =[a_1,\,\ldots,\,a_j]$, by the mirror formula
    \eqref{mirror} we have
    $[a_j,\,\ldots,\,a_1] = \frac {q_{j-1}}{q_j}$ for every
    $1\leq j \leq k$. In particular $\frac rn= \frac {p_k}{q_k}$ and
    $\frac {\tilde r}n= \frac {q_{k-1}}{q_k}$.
    \\[0.15cm]
    By adapting some well known facts about  the slow-additive-Farey continued fraction
    algorithm to the present context, we start setting
    \[
        \Psi_h=\Phi_1^{h-1}\Phi_0=\begin{pmatrix} 0 & 1 \\ 1 &
            h \end{pmatrix}
    \]
    so that thanks to Equation~\eqref{eq:rational2} we can write
    \[
        {p_k\choose q_{k}} = \Psi_{a_1}\cdots \Psi_{a_{k}} {0\choose
          1} \quad \text{and} \quad {q_{k-1}\choose q_{k}} =
        \Psi_{a_k}\cdots \Psi_{a_{1}} {0\choose 1}.
    \]
    Since, as already noted, both $a_1$ and $a_k$ are larger than 1,
    we also have
    \[
        {p_k\choose q_{k}} = \Psi_{a_1}\cdots
        \Psi_{a_{k-1}}\Phi_1^{a_{k}-1} {1\choose 1}
        \quad \text{and} \quad
        {q_{k-1}\choose q_{k}} = \Psi_{a_k}\cdots \Psi_{a_2}
        \Phi_1^{a_{1}-1} {1\choose 1}
    \]
    We now show that the $\dep{r/n}-1=\sum_{j=1}^k a_j -2$ pairs of
    partitions of $n$ with two different parts generated by the dual
    pair $\frac rn$ and $\frac {\tilde r}n$ can be obtained as
    \begin{equation}\label{rappresentazione}
        {p_k\choose q_{k}} = \Psi_{a_1}\cdots \Psi_{a_{k-j-1}}
        \Phi_{1}^{a_{k-j}- r-1}{n_2\choose n_1} \quad \text{and} \quad
        {q_{k-1}\choose q_{k}} = \Psi_{a_k}\cdots \Psi_{a_{k-j+1}}
        \Phi_{1}^{r} {{\tilde n_2}\choose {\tilde n_1}}
    \end{equation}
    for some $j=0,\,\ldots,\,k-1$ and $r=0,\,\ldots,\,a_{k-j}-1$ (with
    $r\geq 1$ for $j=0$ and $r\leq a_1-2$ for $j=k-1$). We first note
    that the choice $(j,r)=(0,0)$ yields the dual pairs
    $(1^n) \vdash n$ and $(n^1) \vdash n$.  For $(j,r)=(0,1)$, a
    straightforward calculation gives
    $$
    {p_k\choose q_{k}} = \Psi_{a_1}\cdots \Psi_{a_{k-1}}
    \Phi_{1}^{a_{k}-2}{1\choose 2}=\begin{pmatrix} p_k -2p_{k-1} &
        p_{k-1} \\ q_k-2q_{k-1} & q_{k-1} \end{pmatrix}{ 1\choose 2}
    $$
    and
    $$
    {q_{k-1}\choose q_k} =  \Phi_{1} {q_{k-1} \choose {q_k - q_{k-1}}}=
    \begin{pmatrix} 1 & 0 \\ 1 & 1 \end{pmatrix} {q_{k-1}\choose {q_k
        - q_{k-1} }}
    $$
    thus producing the pair of dual partitions
    $(2^{q_{k-1}}, 1^{q_k-2q_{k-1}})\vdash n$ and
    $(q_k-q_{k-1}, q_{k-1})\vdash n$. We can now proceed by
    induction. Suppose that, for some choice of $(j,r)$,
    (\ref{rappresentazione}) produces a pair of dual partitions
    $(n_1^{k_1}, n_2^{k_2})$ and
    $({\tilde n}_1^{{\tilde k}_1}, {\tilde n}_2^{{\tilde k}_2})$, that
    is
    \[
        \Psi_{a_1}\cdots \Psi_{a_{k-j-1}} \Phi_{1}^{a_{k-j}- r-1}{n_2\choose n_1}= \begin{pmatrix} t & s \\  k_2 & k_1 \end{pmatrix}{n_2\choose n_1}
    \]
    and
    \[
        \Psi_{a_k}\cdots \Psi_{a_{k-j+1}} \Phi_{1}^{r} {{\tilde
            n_2}\choose {\tilde n_1}} =\begin{pmatrix} {\tilde t }&
            {\tilde s} \\ {\tilde k_2} & {\tilde
              k_1} \end{pmatrix}{{\tilde n_2}\choose {\tilde n_1}}.
    \]
    Assuming $r<a_{k-j}-1$ we can make the transition from $(j,r)$ to
    $(j,r+1)$ and get
    \[
        \Psi_{a_1}\cdots \Psi_{a_{k-j-1}} \Phi_{1}^{a_{k-j}- r-2}{n_2\choose n_1+n_2}= \begin{pmatrix} t-s & s \\  k_2-k_1 & k_1 \end{pmatrix}{n_2\choose n_1+n_2}
    \]
    and
    \[
        \Psi_{a_k}\cdots \Psi_{a_{k-j+1}} \Phi_{1}^{r+1} {{\tilde
            n_2}\choose {\tilde n_1}-{\tilde n_2}} =\begin{pmatrix}
            {\tilde t }+{\tilde s} & {\tilde s} \\ {\tilde k_2}+
            {\tilde k_1} & {\tilde k_1} \end{pmatrix}{{\tilde
            n_2}\choose {\tilde n_1}-{\tilde n_2}} =
        \begin{pmatrix} {\tilde t }+{\tilde s} & {\tilde s} \\ n_1 &
            n_2 \end{pmatrix}{k_1 \choose k_2},
    \]
    where the last identity takes into account the duality of the
    previous pair. As it can be easily checked, we obtain a new pair
    of dual partitions of $n$.  A similar argument applies when
    $r=a_{k-j}-1$ and the transition is from $(j,r)$ to $(j+1,0)$.
\end{proof}

\begin{example}
    We take $n=11$, a prime number, so that we know that all the
    $p(2,11)=27$ partitions into two parts can be generated from the
    dynamics of the Farey map. In the following table we show the
    orbits generated by $\frac r{11}$, where we let
    $r=1,\,\ldots,\,5$, because Lemma~\ref{lemma:paired} implies that
    for $r=6,\,\ldots,\,10$ the orbit generated by $\frac r{11}$ is
    pointwise the same as that generated by $\frac{11-r}{11}$.
    \[
        \begin{array}{c|c|c|c|c|c}
            \hspace{1cm}m\hspace{1cm} & r=1 &r=2  & r=3  & r=4  &r=5 \\
            \hline
            \hline
            1  & (10,1) & (9,2)  & (8,3)  & (7,4)  &(6,5) \\
            \hline
            2  & (9,1^2) & (7,2^2)  & (5,3^2)  & (4^2,3)  &(5^2,1) \\
            \hline
            3  & (8,1^3) & (5,2^3)  & (3^3,2)  & (3^3,1^2)  &(4^2,1^3) \\
            \hline
            4  & (7,1^4) & (3,2^4)  & (2^4,1^3)  & (2^3,1^5)  &(3^2,1^5) \\
            \hline
            5  & (6,1^5) & (2^5,1)  &   &  &(2^2,1^7) \\
            \hline
            6  & (5,1^6) &  &   &  & \\
            \hline
            7  & (4,1^7) &  &   &  & \\
            \hline
            8  & (3,1^8) &  &   &  & \\
            \hline
            9  & (2,1^9) &  &   &  & \\
            \hline
        \end{array}
    \]
    This table should clarify the structure of conjugate
    partitions. For instance, the fractions $\frac 2{11}$ and
    $\frac 5{11}$ generate orbits of dual partitions, in the sense of
    the previous Theorem~\ref{thmconj}. That is, if we read the column
    $r=2$ from top to bottom we find the dual partitions of the column
    $r=5$ read from bottom to top. The same holds for $r=3$ and
    $r=4$. For $r=1$ we have an orbit of partitions which is
    self-dual, meaning that the dual of the $m$-th of the column
    partition is the $(10-m)$-th one in the same column.
\end{example}

\subsection{The extended Farey map} \label{sec:extended-farey}

It will now be more convenient to write partitions as
$(n_1^{k_1}, n_2^{k_2}) \vdash n$ with $n_1>n_2\geq 1$ as
\[
    (n_1, n_2) \times [k_1, k_2]\vdash n.
\]
In part, this is to some extent placing the numbers $n_1$ and $n_2$ on
the same footing as the multiplicities $k_1$ and $k_2$.

We start with extending the definition of the original Farey map $F$
to what we call the \emph{extended Farey map}, which acts on
partitions as follows:
\begin{align*}
    \tilde{F}((n_1,n_2) \times [k_1, k_2]) &=
    \begin{cases}
        \tilde{F}_0((n_1,n_2) \times [k_1, k_2])\, , &\, \text{if } n_2 \geq n_1-n_2\\
        \tilde{F}_1((n_1,n_2) \times [k_1, k_2])\, , &\, \text{if } n_1-n_2 \geq n_2
    \end{cases}\\
    &= \begin{cases}
        (n_2, n_1-n_2) \times  [k_1+k_2, k_1]\, , &\, \text{if } n_2 \geq n_1-n_2\\
        (n_1-n_2 , n_2) \times [k_1, k_1+ k_2]\, , &\, \text{if } n_1-n_2 \geq n_2
    \end{cases}
\end{align*}
The action of $\tilde F$ on $(n_1,n_2)$ is just the action of the
Farey map written as in Section~\ref{sec:Fmatrix}. The corresponding
action on the multiplicities $(k_1,k_2)$ is the one obtained as
follows. We can write
\[
    n = k_1n_1+k_2n_2 = \vvec{k_2}{k_1}^\top\vvec{n_2}{n_1},
\]
so that if the action on $(n_1,n_2)$ is given by $F_i$, with $i=0$ or
$1$, then
\[
    n = \vvec{k_2}{k_1}^\top\Phi_iF_i\vvec{n_2}{n_1} =
    \left(\Phi_i^\top\vvec{k_2}{k_1}\right)^\top\left(F_i\vvec{n_2}{n_1}\right)
\]
and the action on the multiplicities is that of $\Phi_i^\top$.

\begin{remark}
    In our work for this paper we quickly arrived at the above map
    $\tilde F$, by simply finding out how the multiplicities should
    transform. But once we wrote it down, we saw that we were
    simply reproducing the natural extension of the Farey map, as
    described, for example, by Arnoux and Nogueira
    \cite{Arnoux-Nogueira-93}. Natural extensions are a standard tool
    in dynamical systems which  change $n$ to $1$ maps into $1$ to $1$
    maps by  extending the dimension of the domain. This gives as an
    interpretation of the extened Farey map $\tilde{F}$ in terms of
    matrix multiplication if we write the partition
    $(n_1,n_2) \times [k_1, k_2]$ as the four-dimensional vector
    $(n_2,n_1,k_2,k_1)$. Indeed we have
    \begin{align*}
        \tilde{F} \begin{pmatrix} n_2 \\ n_1 \\ k_2 \\ k_1 \end{pmatrix}
        &=
        \begin{cases}
            \begin{pmatrix} F_0 & 0\\0& \Phi_0^\top\end{pmatrix}
            \begin{pmatrix} n_2 \\ n_1 \\ k_2 \\ k_1 \end{pmatrix}\, ,
            &
            \, \text{if } n_2\geq n_1-n_2\\[1cm]
            \begin{pmatrix} F_1 & 0\\0& \Phi_1^\top\end{pmatrix}
            \begin{pmatrix} n_2 \\ n_1 \\ k_2 \\ k_1 \end{pmatrix}\, ,
            &
            \, \text{if } n_1-n_2\geq n_2\\
        \end{cases}\\
        &= \begin{cases}
            \begin{pmatrix}
                -1&1&0&0\\1&0&0&0 \\ 0&0& 0&1\\ 0&0&1&1
            \end{pmatrix}
            \begin{pmatrix} n_2 \\ n_1 \\ k_2 \\ k_1 \end{pmatrix}\, , &
            \, \text{if } n_2\geq n_1-n_2\\[1cm]
            \begin{pmatrix} 1&0& 0 & 0\\-1&1&0&0 \\ 0&0& 1&1\\
                0&0&0&1\end{pmatrix}
            \begin{pmatrix} n_2 \\ n_1 \\ k_2 \\ k_1 \end{pmatrix}\, ,
            &
            \, \text{if } n_1-n_2\geq n_2\\
        \end{cases}
        =
        \begin{cases}
            \begin{pmatrix} n_1-n_2 \\ n_2 \\ k_1 \\ k_1+k_2 \end{pmatrix}\, ,
            &
            \, \text{if } n_2\geq n_1-n_2\\[1cm]
            \begin{pmatrix} n_2 \\ n_1-n_2 \\ k_1+k_2 \\ k_1 \end{pmatrix}\, ,
            &
            \, \text{if } n_1-n_2\geq n_2\\
        \end{cases}
    \end{align*}
    We wrote this fully out as this will directly generalize to higher
    dimensions, as we will see in Section~\ref{npartitions}.
\end{remark}

The extended Farey map $\tilde{F}$ maps a partition of $n$ to a new
partition of $n$. Indeed if $(n_1, n_2) \times [k_1, k_2]\vdash n$, a
simple calculation shows that both
\[
    (n_2, n_1-n_2) \times [k_1+k_2, k_1] \quad\text{and}\quad (n_1-n_2
    , n_2) \times [k_1, k_1+ k_2]
\]
are again partitons of $n$. An immediate consequence of the definition
given in Section~\ref{2-intro} and of the present construction is that
if $r$ is an integer such that $1\leq r<n$ and $(r,n) = 1$ then
repeatedly applying $\tilde F$ to $(n,r)\times[1,0]$ yields exactly
the orbit of partitions generated by dynamics of the Farey map
starting from $\frac rn$. In particular, the partition
$(n,r)\times[1,0]$ eventually maps to $(1,1)\times [k_1,k_2]$ for some
$k_1$ and $k_2$.

\begin{example}
    \label{ex8/19i}
    In Example~\ref{ex8/19} we considered the orbit of partitions of
    $19$ generated by $\frac 8{19}$. We now show the same orbit of
    partitions as it can be obtained through the map $\tilde F$:
    \begin{align*}
        (19, 8) \times [ 1 , 0 ] &\overset{\tilde{F}_1}{\longmapsto}%
        (11, 8) \times [ 1 , 1 ] \overset{\tilde{F}_0}{\longmapsto}%
        (8, 3) \times [ 2 , 1 ] \overset{\tilde{F}_1}{\longmapsto}%
        (5,3) \times [ 2 , 3 ] \\
        &\overset{\tilde{F}_0}{\longmapsto}%
        (3, 2) \times [ 5 , 2 ] \overset{\tilde{F}_0}{\longmapsto}%
        (2, 1) \times [ 7 , 5 ] \overset{\tilde{F}_0}{\longmapsto}%
        (1 ,1) \times [12 , 7].
    \end{align*}
    In the last step, exactly as in the construction of
    Section~\ref{2-intro}, we could have used $\tilde{F}_1$ obtaining
    the partition $(1,1)\times [7,12]$.
\end{example}

In the more general setting of the present section there is no need to
restrict ourselves to starting with a vector $(n,r)$ of relatively
prime numbers, nor with the vector $[1,0]$ for the
multiplicities. This will allow us to consider also partitions with
non-coprime numbers and/or multiplicities and will be key to proving
the formula for $p(2,n)$ in Section~\ref{p(2,n)}.

\begin{lemma}
    \label{lemma:de}
    Let $n\geq 2$ and $1\leq r<n$. It holds that
    $\tilde{F}^m((n,r)\times [1,0]) = (n_1,n_2)\times [k_1,k_2]$ if
    and only if
    $\tilde{F}^m((dn,dr)\times [e,0]) = (dn_1,dn_2)\times [ek_1,ek_2]$
    for every $d\geq 1$ and $e\geq 1$.
\end{lemma}
\begin{proof}
    The map $\tilde F$ is linear in all of its arguments.
\end{proof}

\noindent As a consequence of the previous lemma, we can now show that
the iterations of the extended Farey map stop at some point. Suppose
that $d=(n,r)$, let $n'=\frac nd$ and $r'=\frac rd$, so that
$(n',r')=1$. We know that
\[
    \tilde{F}^{m}((n',r')\times [1,0]) = (1,1)\times [k_1,k_2]
\]
for some $k_1$ and $k_2$, and $m=\dep{r'/n'}=\dep{r/n}$. Thus from
Lemma~\ref{lemma:de} it follows that
\[
    \tilde{F}^m((n,r)\times [e,0]) = (d,d)\times [ek_1,ek_2].
\]
In other words, the starting partition $(n,r)\times [e,0]$ is
eventually mapped to a partition having equal numbers.

\begin{definition}
    We shall call the sequence of partitions $F^m((n,r)\times [e,0])$
    with $m=1,\,\ldots,\,\dep{\frac rn}$ the \emph{orbit of partitions
      generated by $(n,r)\times [e,0]$} under the extended Farey map.
\end{definition}

\begin{definition}
    Suppose that $(n_1, n_2) \times [k_1,k_2]$ eventually maps to
    $(m_1, m_2) \times [l_1,l_2]$ under the extended Farey map.  Then
    we say that $(n_1, n_2) \times [k_1,k_2]$ is an \emph{ancestor} of
    $(m_1, m_2) \times [l_1,l_2]$ and that
    $(m_1, m_2) \times [l_1,l_2]$ is a \emph{descendant} of
    $(n_1, n_2) \times [k_1,k_2]$.
\end{definition}

We now show how Lemma~\ref{lemma:de} works in practice with some
examples.

\begin{example}
    We saw in Example~\ref{ex8/19i} that the partition
    $ (5, 3) \times [ 2 , 3 ] \vdash 19$ appears in the third generation
    starting with $(19, 8) \times [ 1 , 0 ]$ at generation
    $0$.
    \\[0.15cm]
    (i) We first consider the partition
    \[
        (10, 6) \times [ 2 , 3 ] \vdash 38,
    \]
    which is simply
    $(2\cdot 5, 2\cdot 3) \times [ 2 , 3 ] \vdash 2\cdot 19$. We now
    show that it appears in the third generation in the orbit of partitions
    of $(38,16)\times [1,0]$. Indeed we have
    \begin{align*}
        (38, 16) \times [ 1 , 0 ] &\overset{\tilde{F}_1}{\longmapsto}%
        (22, 16) \times [ 1 , 1 ] \overset{\tilde{F}_0}{\longmapsto}%
        (16, 6) \times [ 2 , 1 ] \overset{\tilde{F}_1}{\longmapsto}%
        (10, 6) \times [ 2 , 3 ]\\
        &\overset{\tilde{F}_0}{\longmapsto}%
        (6, 4) \times [ 5 , 2 ] \overset{\tilde{F}_0}{\longmapsto}%
        (4, 2) \times [ 7 , 5 ] \overset{\tilde{F}_0}{\longmapsto}%
        (2 ,2) \times [12 , 7],
    \end{align*}
    which is the same orbit as for $(19, 8) \times [ 1 , 0 ]$ save for
    multiplying all the numbers by $2$.
    \\[0.15cm]
    (ii) Now we consider
    \[
        (5, 3) \times [ 6 , 9 ] \vdash 57,
    \]
    which is
    $ (5, 3) \times [ 3\cdot 2 , 3\cdot 3 ] \vdash 3 \cdot 19$ and
    appears in the third generation starting from $(19,8)\times[3,0]$:
    \begin{align*}
        (19, 8) \times [ 3 , 0 ] &\overset{\tilde{F}_1}{\longmapsto}%
        (11, 8) \times [ 3 , 3 ] \overset{\tilde{F}_0}{\longmapsto}%
        (8, 3) \times [ 6 , 3 ] \overset{\tilde{F}_1}{\longmapsto}%
        (5,3) \times [ 6 , 9 ] \\
        &\overset{\tilde{F}_0}{\longmapsto}%
        (3, 2) \times [ 15 , 6 ] \overset{\tilde{F}_0}{\longmapsto}%
        (2, 1) \times [ 21 , 15 ] \overset{\tilde{F}_0}{\longmapsto}%
        (1 ,1) \times [36 , 21].
    \end{align*}
    This orbit is the same as for $(19, 8) \times [ 1 , 0 ]$, save for
    now multiplying all the multiplicities by $3$.
    \\[0.15cm]
    (iii) Finally, consider
    \[
        (10, 6) \times [ 6 , 9 ] \vdash 114,
    \]
    which is
    $( 2\cdot 5, 2\cdot 3) \times [ 3\cdot 2 , 3\cdot 3 ] \vdash 2\cdot
    3 \cdot 19$. We have
    \begin{align*}
        (38, 16) \times [ 3 , 0 ] &\overset{\tilde{F}_1}{\longmapsto}%
        (22, 16) \times [ 3 , 3 ] \overset{\tilde{F}_0}{\longmapsto}%
        (16, 6) \times [ 6 , 3 ] \overset{\tilde{F}_1}{\longmapsto}%
        (10,6) \times [ 6 , 9 ] \\
        &\overset{\tilde{F}_0}{\longmapsto}%
        (6, 4) \times [ 15 , 6 ] \overset{\tilde{F}_0}{\longmapsto}%
        (4, 2) \times [ 21 , 15 ] \overset{\tilde{F}_0}{\longmapsto}%
        (2 ,2) \times [36 , 21],
    \end{align*}
    ending up with the same orbit as $(19, 8) \times [ 1 , 0 ]$, but
    now multiplying all the numbers by $2$ and all the multiplicities
    by $3$.
\end{example}

For orbits of partitions generated with the extended Farey map we have
this more general version of Lemma~\ref{lemma:paired}.

\begin{prop}
    \label{uniqueness}
    The partitions $(n,r) \times [e,0]$ and $(n,n-r) \times [e,0]$
    generate the same orbit of partitions.
\end{prop}

\begin{proof}
    If $r= n-r$ then the result is obvious. So we can assume that
    $r>n-r$.  Thus we apply $\tilde{F}_0$ to $(n,r) \times [e,0]$ and
    $\tilde{F}_1$ to $(n,n-r) \times [e,0]$ getting
    \begin{align*}
        \tilde{F}_0((n,r) \times [e,0] ) &= (r, n-r) \times [e,e] = \\
        &= (n- (n-r), n-r) \times [e,e] = \tilde{F}_1((n,n-r) \times
        [e,0]).
    \end{align*}
    Thus after the first application of the extended Farey map, we
    have the same partition.
\end{proof}

\subsection{Conjugation and palindromes, again}
\label{palindrome 2 dim}

The extended Farey map reflects and respects conjugation.

\begin{prop}
    \label{prop:conj-maps}
    The diagram
    \begin{eqnarray*}
        (n_1,n_2) \times [k_1, k_2] & \sim_{\mathcal{C}} & (k_1+k_2,
        k_1 ) \times [n_2, n_1-n_2]
        \\
        \tilde{F}_0 \downmapsto & &\upmapsto\tilde{F}_0
        \\
        (n_2, n_1-n_2) \times [k_1+k_2, k_1] & \sim_{\mathcal{C}} &
        (2k_1+k_2, k_1+k_2) \times [n_1-n_2, 2n_2-n_1]
    \end{eqnarray*}
    when $n_2 \geq n_1 - n_2$, and the diagram
    \begin{eqnarray*}
        (n_1,n_2) \times [k_1, k_2] & \sim_{\mathcal{C}} & (k_1+k_2,
        k_1 ) \times [n_2, n_1-n_2]
        \\
        \tilde{F}_1 \downmapsto & &\upmapsto\tilde{F}_1        \\
        (n_1-n_2, n_2) \times [k_1, k_1+ k_2 ] & \sim_{\mathcal{C}} &
        (2k_1+k_2, k_1) \times [n_2, n_1-2n_2]
    \end{eqnarray*}
    when $n_2\leq n_1-n_2$, are both commutative.
\end{prop}
\begin{proof}
    It is a simple verification using the definition of the map
    $\tilde F$ and the conjugation rule.
\end{proof}

\noindent Commutative diagrams as the ones shown in the above
proposition can be glued together following an orbit of partitions. As
an example, consider
\begin{eqnarray*}
    (19, 15) \times [1,0]  &  \sim_{\mathcal{C}}    &   (1,1) \times [15, 4] \\
    \tilde{F}_0 \downmapsto & &\upmapsto\tilde{F}_0\\
    (15, 4) \times [1,1]  &  \sim_{\mathcal{C}}    &   (2,1) \times [4, 11] \\
    \tilde{F}_1 \downmapsto & &\upmapsto\tilde{F}_1        \\
    (11, 4) \times [1,2]  &  \sim_{\mathcal{C}}    &   (3,1) \times [4, 7] \\
    \tilde{F}_1 \downmapsto & &\upmapsto\tilde{F}_1        \\
    (7, 4) \times [1,3]  &  \sim_{\mathcal{C}}    &   (4,1) \times [4, 3] \\
    \tilde{F}_0 \downmapsto & &\upmapsto\tilde{F}_0\\
    (4, 3) \times [4,1]  &  \sim_{\mathcal{C}}    &   (5,4) \times [3, 1] \\
    \tilde{F}_0 \downmapsto & &\upmapsto\tilde{F}_0\\
    (3,1) \times [5,4]  &  \sim_{\mathcal{C}}    &   (9,5) \times [1, 2] \\
    \tilde{F}_1 \downmapsto & &\upmapsto\tilde{F}_1        \\
    (2, 1) \times [5,9]  &  \sim_{\mathcal{C}}    &   (14,5) \times [1, 1] \\
    \tilde{F}_0 \downmapsto & &\upmapsto\tilde{F}_0\\
    (1, 1) \times [14,5]  &  \sim_{\mathcal{C}}    &   (19,14) \times [1, 0] \\
\end{eqnarray*}
We remark that the last two lines could also have been
\begin{eqnarray*}
    (2, 1) \times [5,9]  &  \sim_{\mathcal{C}}    &   (14,5) \times [1, 1] \\
    \tilde{F}_1 \downmapsto & &\upmapsto\tilde{F}_1        \\
    (1, 1) \times [5,14]  &  \sim_{\mathcal{C}}    &   (19,5) \times [1, 0].
\end{eqnarray*}
Hence the act of conjugation can be viewed as simply reversing the
arrows of the extended map, and this can be described by another
version of Theorem~\ref{thmconj}.

\begin{theorem}[Palidromes Version 2]
    \label{2-palindorme}
    Suppose that
    $ \tilde{F}_{\sigma_1} , \ldots , \tilde{F}_{\sigma_\ell}$ is a
    sequence of extended Farey maps such that
    \[
        \tilde{F}_{\sigma_\ell}\circ \cdots \circ
        \tilde{F}_{\sigma_1}((n_1, n_2 ) \times [k_1, k_2] )=
        (\bar{n}_1, \bar{n}_2) \times [\bar{k}_1 , \bar{k}_2].
    \]
    Then
    \[
        \tilde{F}_{\sigma_1}\circ \cdots \circ
        \tilde{F}_{\sigma_\ell}( (\bar{k}_1 + \bar{k}_2, \bar{k}_1)
        \times [\bar{n}_2 , \bar{n}_1-\bar{n}_2] )=(k_1+k_2, k_1 )
        \times [n_2, n_1-n_2].
    \]
\end{theorem}
\begin{proof}
    It is a repeated application of
    Proposition~\ref{prop:conj-maps}. To shorten the notation, for
    $j=1,\,\ldots,\,\ell$ denote by $\lambda^{(j)}$ the partition
    $\tilde{F}_{\sigma_j}\circ \cdots \circ \tilde{F}_{\sigma_1}((n_1,
    n_2 ) \times [k_1, k_2] )$ and by $\lambda^{(j)}_{\mathcal C}$ its
    conjugate. Now let $1\leq j< \ell$ and suppose that
    \[
        \tilde{F}_{\sigma_1}\circ \cdots \circ
        \tilde{F}_{\sigma_{j-1}}(\lambda^{(j)}_{\mathcal C}) =
        \lambda^{(1)}_{\mathcal{C}}.
    \]
    Since $\lambda^{(j+1)} = \tilde{F}_{\sigma_j}(\lambda^{(j)})$,
    Proposition~\ref{prop:conj-maps} yields
    $\lambda^{(j)}_{\mathcal C} =
    \tilde{F}_{\sigma_j}(\lambda^{(j+1)}_{\mathcal C})$, so that
    \[
        \lambda^{(1)}_{\mathcal{C}} = \tilde{F}_{\sigma_1}\circ \cdots
        \circ \tilde{F}_{\sigma_{j-1}}(\lambda^{(j)}_{\mathcal C}) =
        \tilde{F}_{\sigma_1}\circ \cdots \circ
        \tilde{F}_{\sigma_{j}}(\lambda^{(j+1)}_{\mathcal C}),
    \]
    proving the inductive step.
\end{proof}

As an application of the previous result about conjugation, we give a
new version of Theorem~\ref{2parts-farey}.

\begin{theorem}
    Suppose that $(n_1,n_2) \times [k_1,k_2]\vdash n$ with
    $(n_1, n_2)=1$ and $(k_1, k_2)=1$. Then there exists some positive
    integer $r$ relatively prime to $n$ and such that
    $(n,r) \times [1,0]$ is an ancestor of
    $(n_1,n_2) \times [k_1,k_2]$.
\end{theorem}

\begin{proof}
    We start with the partition $(n_1,n_2) \times [k_1,k_2]\vdash n$
    and consider its conjugate partition, namely
    $(k _1 + k_2, k_1) \times [n_2, n_1-n_2]$, which is still a
    partition of $n$. Since $(k_1+k_2, k_1)=1$ we know that the
    fraction $\frac{k_1}{k_1+k_2}$ appears in some level of the Farey
    tree, thus there is a sequence of Farey matrices
    $F_{\sigma_1},\,\ldots,\,F_{\sigma_\ell}$ such that
    \[
        F_{\sigma_{\ell}}\cdots F_{\sigma_1}\vvec{k_1}{k_1+k_2} =
        \vvec{1}{1}.
    \]
    Hence for the extended maps we have
    \[
        \tilde{F}_{\sigma_\ell}\circ\cdots\circ \tilde{F}_{\sigma_1}
        ((k _1 + k_2, k_1) \times [n_2, n_1-n_2]) = (1,1) \times [r,s]
    \]
    with $(1,1) \times [r,s]\vdash n$. Note that $r+s=n$ and that
    Lemma~\ref{lemma:de} implies that $(r,n)=1$. Now, the partition
    $(1,1) \times [r,s] $ is conjugate to
    $(r+s,r) \times [1,0]= (n,r) \times [1,0]$, thus by the Power of
    Palindromes we have
    \[
        \tilde{F}_{\sigma_1} \cdots \tilde{F}_{\sigma_\ell} ( (n,r)
        \times [1, 0]) = ( n_1 , n_2) \times
        [k_1,k_2]
    \]
    and we are done.
\end{proof}

\subsection{Formula for $p(2,n)$}
\label{p(2,n)}

We now use the previous results to describe the partitions which can be
obtained through the extended Farey map and derive the formula for
$p(2,n)$.

\begin{theorem}\label{two sum}
    Let $n\geq 2$ be an integer. Every partition of $n$ can be
    obtained from the dynamics of the extended Farey map $\tilde
    F$.
\end{theorem}
\begin{proof}
    Consider a partition $(n_1,n_2)\times [k_1,k_2]$ of $n$ and let
    $d=(n_1,n_2)$ and $e = (k_1,k_2)$. By setting
    $n_1'=\frac{n_1}{d}$, $n_2'=\frac{n_2}{d}$, $k_1'=\frac{k_1}{e}$,
    $k_2'=\frac{k_2}{e}$ we have that
    \[
        (n_1',n_2')\times [k_1',k_2'] \vdash \frac{n}{de}
    \]
    is by construction a partition with relatively prime numbers and
    multiplicities. Thus there exists $h$ such that
    $\left(\frac n{de},h\right)\times [1,0]$ is an ancestor of
    $(n_1',n_2')\times [k_1',k_2']$ . Hence by Lemma~\ref{lemma:de} we
    have that
    \[
        \left(\frac ne,hd\right)\times [e,0]
    \]
    is an ancestor of
    $(dn_1',dn_2')\times [ek_1',ek_2']=(n_1,n_2)\times [k_1,k_2]$.
\end{proof}

\begin{theorem}\label{p2nformula}
    Let $n\geq 2$ be an integer number. Then
    \[
        p(2,n) = \frac{1}{2}
        \sum_{r=1}^{n-1} \left(\dep{\frac rn} -1 \right) \sigma_0
        ((r,n)).
    \]
\end{theorem}
\begin{proof}
    To shorten the notation, we denote by $\mathcal{O}(\lambda)$ the
    orbit of partition generated by the extended Farey map starting
    from the partition $\lambda$. We claim that the set of the
    partitions of $n$ into two parts is
    \[
        P = \bigcup_{1\leq r\leq \frac n2} \bigcup_{e\mid (r,n)}
        \mathcal{O}\left( \left(\frac ne,\frac re\right)\times
          [e,0]\right).
    \]
    It is clear that each partition of this set is a partition of $n$
    into two parts. For the converse we consider a partition
    $(n_1,n_2)\times [k_1,k_2]$ of $n$. By arguing as in the proof of
    Theorem~\ref{two sum} we have that there exists $h$ such that
    $\left(\frac ne,hd\right)\times [e,0]$ is an ancestor of our
    partition. Furthermore, it is possible to choose $h$ such that
    $1\leq h\leq \frac{n}{2de}$. Thus our partition
    $(n_1,n_2)\times [k_1,k_2]$ is in the orbit of
    $\left(\frac ne,\frac re\right)\times [e,0]$, where $r=hde$. We
    have that $de\leq r\leq \frac n2$ and also that $e$ is a divisor
    of $(r,n)$ since it divides both $n$ and $r$, thus our partition
    is in the set $P$.
    \\[0.15cm]
    To prove the formula for $p(2,n)$ it now suffices to count the
    elements of the above set $P$. For each $1\leq r\leq \frac n2$ we
    consider the $\sigma_0((r,n))$ disjoint orbits generated by
    $\left(\frac ne,\frac re\right)\times [e,0]$. Each of them
    contains
    \[
        \dep{\frac{r/e}{n/e}}-1 = \dep{\frac rn}-1
    \]
    pairwise distinct partitions of $n$. Thus the number of partitions
    in $P$ is
    \[
        \sum_{1\leq r\leq \frac n2} \sum_{e\mid (r,n)}
        \left(\dep{\frac rn}-1\right) =\sum_{1\leq r\leq \frac n2}
        \left(\dep{\frac rn}-1\right) \sigma_0((r,n)).
    \]
    If we extend the first sum over $1\leq r\leq n-1$ we are counting
    each partition twice due to Proposition~\ref{uniqueness}, for
    which $(n,r) \times [e,0] $ give rise to the same descendants as
    $(n,n-r) \times [e,0]$.
\end{proof}

\begin{example}
    Here are all the partitions of $n=12$ into two numbers:
    \[
        \begin{array}{cccc}
            (11,1 ) \times [1,1]  & (10,2 ) \times [1,1] & (10,1 ) \times [1,2] & (9,3 ) \times [1,1] \\

            (9,1 ) \times [1,3]  & (8,4 ) \times [1,1] & (8,2 ) \times [1,2] & (8, 1 ) \times [1,4] \\
            (7,5 ) \times [1,1]  & (7,1 ) \times [1,5] & (6,3 ) \times [1,2] & (6,2 ) \times [1,3] \\
            (6,1 ) \times [1,6]  & (5,2)\times [2,1] & (5,1) \times [2, 2] & (5,1 ) \times [1,7] \\
            (4, 2 ) \times [2,2] & (4,1 ) \times [2,4]  &(4,2 ) \times [1,4]  & (4,1 ) \times [1,8] \\
            (3,1 ) \times [3,3] & (3,2 ) \times [2,3] & (3,1 ) \times [2,6]  & (3, 1 ) \times [1,9] \\
            (2,1 ) \times [5,2] & (2,1 ) \times [4,4] & (2,1 ) \times [3,6]  & (2,1 ) \times [2,8] \\
            (2, 1 ) \times [1,10] & \\
        \end{array}
    \]
    Thus $p(2,12) = 29$. Indeed we have
    \begin{align*}
        29 &= \left( \dep{ \frac{1}{12} } -1 \right) \sigma_0
        (\gcd (1,12)) +
        \left( \dep{ \frac{2}{12} } -1 \right) \sigma_0 (\gcd (2,12))+~ \\
        &\hspace{0.7cm} + \left( \dep{ \frac{3}{12} } -1 \right) \sigma_0
        (\gcd (3,12))
        +  \left( \dep{ \frac{4}{12} } -1 \right) \sigma_0 (\gcd (4,12))+~\\
        &\hspace{1.1cm} + \left( \dep{ \frac{5}{12} } -1 \right) \sigma_0
        (\gcd (5,12))
        +  \left( \dep{ \frac{6}{12} } -1 \right) \sigma_0 (\gcd (6,12)).
    \end{align*}
    Let us see how each fraction $\frac r{12}$ with $1\leq r\leq 6$
    will produce
    $\left( \dep{ \frac{r}{12} } -1 \right) \sigma_0 ((r,12))$
    distinct partitions of $12$ into two parts.  We start with the two
    values of $r$ for which $\sigma_0 ((r,12)) = 1$ namely $r=1$ and
    $r=5$. We have
    \begin{align*}
        (12, 1)\times [1,0] &\overset{ \tilde{F}_1}{\longmapsto} (11,
        1)\times [1,1] \overset{ \tilde{F}_1}{\longmapsto} (10,
        1)\times [1,2]
        \\
        &\overset{ \tilde{F}_1}{\longmapsto} (9, 1)\times [1,3]
        \overset{ \tilde{F}_1}{\longmapsto} (8, 1)\times [1,4]
        \overset{ \tilde{F}_1}{\longmapsto} (7, 1)\times [1,5]
        \\
        &\overset{ \tilde{F}_1}{\longmapsto} (6, 1)\times [1,6]
        \overset{ \tilde{F}_1}{\longmapsto} (5, 1)\times [1,7]
        \overset{ \tilde{F}_1}{\longmapsto} (4, 1)\times [1,8]
        \\
        &\overset{ \tilde{F}_1}{\longmapsto} (3, 1)\times [1,9]
        \overset{ \tilde{F}_1}{\longmapsto} (2, 1)\times [1,10]
    \end{align*}
    This accounts for precisely $\dep{1/12}-1=10$ of the desired
    partitions. Further, each of these will occur uniquely. Similarly,
    for $r=5$ we have
    \[
        (12, 5)\times [1,0] \overset{\tilde{F}_1}{\longmapsto} (7,
        5)\times [1,1] \overset{ \tilde{F}_0}{\longmapsto} (5,
        2)\times [2,1] \overset{ \tilde{F}_1}{\longmapsto} (3,
        2)\times [2,3] \overset{ \tilde{F}_0}{\longmapsto} (2,
        1)\times [5,2]
    \]
    giving us $\dep{5/12} -1=4$ additional partitions. Now we consider
    $r=2$ and see the partitions linked to the fraction $2/12$. In
    this case $(r,n)=2$, so that we have two different choices for
    $e$, namely $e=1$ and $e=2$. Choosing $e=1$ we start with
    \[
        (12, 2)\times [1,0] \overset{ \tilde{F}_1}{\longmapsto}  (10, 2)\times
        [1,1] \overset{ \tilde{F}_1}{\longmapsto}  (8, 2)\times [1,2]
        \overset{ \tilde{F}_1}{\longmapsto}  (6, 2)\times [1,3] \overset{
          \tilde{F}_1}{\longmapsto}  (4,2)\times [1,4]
    \]
    giving us $\dep{2/12} -1 = 4$ additional partitions. Choosing
    $e=2$ we also have
    \[
        (6, 1)\times [2,0] \overset{ \tilde{F}_1}{\longmapsto}  (5, 1)\times
        [2,2] \overset{ \tilde{F}_1}{\longmapsto}  (4, 1)\times [2,4]
        \overset{ \tilde{F}_1}{\longmapsto}  (3, 1)\times [2,6] \overset{
          \tilde{F}_1}{\longmapsto}  (2,1)\times [2,8]
    \]
    giving us another $4$ partitions. Now we consider $r=3$, that is
    the fraction $3/12$. For $e=1$ we have
    \[
        (12, 3)\times [1,0] \overset{ \tilde{F}_1}{\longmapsto}  (9, 3)\times
        [1,1] \overset{ \tilde{F}_1}{\longmapsto}  (6,3)\times [1,2]
    \]
    giving us $\dep{3/12} -1=2$ partitions.  But we also have the
    choice $e=3$, yielding
    \[
        (4, 1)\times [3,0] \overset{ \tilde{F}_1}{\longmapsto}  (3, 1)\times
        [3,3] \overset{ \tilde{F}_1}{\longmapsto}  (2, 1)\times [3,6]
    \]
    giving us $2$ more partitions. The remaining cases $r=4$, $r=5$
    and $r=6$ work in the same way and complete the list of partitions
    of $12$.
\end{example}


\section{On partitions into many parts}
\label{npartitions}
In this section we begin the extension of the construction explained in Sections \ref{section-farey-map} and \ref{2 case}. The Farey map may be defined to act on a cone in $\R^2$ as in Section \ref{sec:Fmatrix}, and we used it to generate partitions into two parts. To generate partitions into $n$  different parts it is necessary to consider a map acting on a subset of $\R^N$. To this aim we consider the Triangle map and its slow version studied in \cite{Garrity1, Bonanno- Del Vigna-Munday}.

\subsection{Background on the additive-slow-Triangle map}
The Triangle map has been introduced in \cite{Garrity1} to define a
type of multidimensional continued fraction
algorithm. Multidimensional continued fractions have been developed
and studied over the years for many reasons (for a background see
Schweiger \cite{Schweiger4} or Karpenkov \cite{Karpenkov13}). Historically, the first
such algorithm (now called the Jacobi-Perron algorithm) was created to
answer a question of Hermite, which was to find an analogue of the
classical fact that a real number has an eventually periodic continued
fraction expansion if and only if the number is a quadratic
irrational. This problem is still open.  Another motivation was to
find methods for good simultaneous Diophantine approximations of
$n$-tuples of real numbers (for example, see Lagarias \cite{Lagarias93}). By now, multidimensional continued
fractions provide a rich source of examples in dynamical systems,
automata theory and many other areas. For background on and
applications of the Triangle map for multidimensional continued
fractions, see \cite{GarrityT05,Garrity1, Karpenkov13, SchweigerF08,
  Schweiger05, Bonanno- Del Vigna-Munday, Fougeron-Skripchenko,
  Bonanno-Del Vigna, Berthe-Steiner-Thuswaldner, Ito}.  Overwhelmingly these papers are concerned with the three dimensional case.

Set
\begin{eqnarray*}
\triangle &:=& \{ (x_1, \ldots, x_n) \in \R^n: 1 >x_1>\cdots > x_n>0\} \\
\triangle_0 &:=& \{ x_1, \ldots, x_n) \in \triangle: x_1+x_n >1\} \\
\triangle_1 &:=& \{ x_1, \ldots, x_n) \in \triangle: x_1+x_n >1\}
\end{eqnarray*}
When $n=2$, we have 

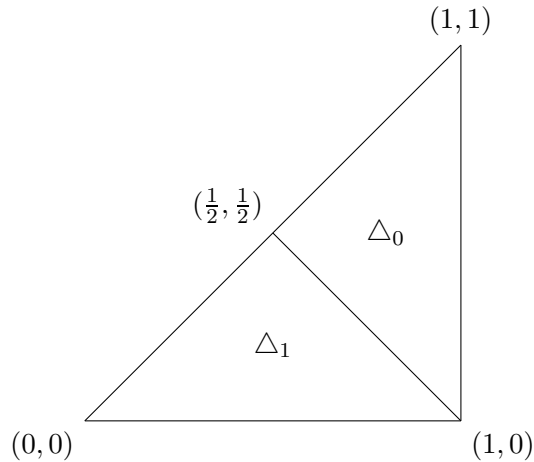
\begin{center}
\begin{tikzpicture}[scale=5]
\draw(0,0)--(1,0);
\draw(0,0)--(1,1);
\draw(1,0)--(1,1);
\draw(1,0)--(1/2,1/2);
\node[] at (4/5,1/2){$\triangle_0$};
\node[] at (1/2,1/5){$\triangle_1$};
\node[below left]at(0,0){$(0,0)$};
\node[below right]at(1,0){$(1,0)$};
\node[above]at(1,1){$(1,1)$};
\node[above left] at (1/2,1/2){$(\frac 12,\frac 12)$};
\end{tikzpicture}
\end{center}

The slow-Triangle map
$T:\triangle_0 \cup \triangle_1 \rightarrow \triangle $ is
\begin{eqnarray*}
T(x_1, \ldots, x_n) &=&  \left\{ \begin{array}{cc} T_0(x_1, \ldots, x_n), &\, \text{if } x_1+x_n>1 \\
                                                       T_1(x_1, \ldots, x_n), &\, \text{if } x_1+x_n<1 \end{array} \right.  \\
                                                       &=&  \left\{ \begin{array}{ccc} \left( \frac{x_2}{x_1},\ldots , \frac{x_n}{x_1}, \frac{1-x_1}{x_1}  \right), &\, \text{if } x_1+x_n>1 \\
                                                        \left( \frac{x_1}{1-x_n}, \ldots, \frac{x_n}{1-x_n}  \right), &\, \text{if } x_1+x_n<1 \end{array} \right.  \\
\end{eqnarray*}
It can be checked that $T_i:\triangle_i \rightarrow \triangle$ is
one-to-one and onto. In analogy with the construction in Section
\ref{farey}, any point $\bar{x} \in \triangle$ is associated to a binary
sequence of zeros and ones $i(\bar{x}) = (i_0,\, i_1,\, i_2,\, \,\ldots)$
by encoding whether its iterations fall in $\triangle_0$ or
$\triangle_1$, that is by the rule $T^{n}(\bar{x}) \in
\triangle_{i_n}$. We call $i(\bar{x})$ the additive-slow-Triangle sequence
of $\bar{x}$. If we concatenate the $1$'s we can associate $\bar{x}$ to a
sequence of nonnegative integers, a sequence that is called the
multiplicative-fast-Triangle sequence and is the analogue of the
continued fraction expansion of a real number. Either of these
sequences tells us a lot about the point $\bar{x}$. For example, when $n=3$,  if the
sequence is eventually periodic, then both $x_1$ and $x_2$ are no worse
than cubic irrationals, both in the same number field of degree less
than or equal to three.

\begin{remark}
    With respect to \cite{Garrity1, Bonanno- Del Vigna-Munday} we have
    not defined the map $T$ on the boundary of $\triangle$ and on the
    hyperplane $x_1+ x_n=1$.   In earlier work, such points are a set of measure zero and hence are ignored.  This creates a problem, though, when we start to link these maps with partitions, as we will discuss in Section \ref{$T_D$}.
\end{remark}

It is natural, and in analogy to Section \ref{sec:Fmatrix} for the
Farey map, to pass from points $(x_1, \ldots ,x_n) $ in $\R^n$ to vectors $(x_0, \ldots, x_n)$
in $\R^{n+1}$ (or points $(x_0, \ldots, x_n)$ in $\R\P^{n+1}$) via sending $(x_1, \ldots ,x_n)  $ to
$(1,x_1, \ldots ,x_n) $ with inverse map $(x_0, \ldots ,x_n) \rightarrow (x_1/x_0, \ldots, x_n/x_0).$ Then, by
an abuse of notation, we set
\begin{eqnarray*}
\triangle &:=& \{(x_0, \ldots, x_n)\in \R^{n+1}: x_0 >x_1>\cdots > x_n>0\} \\
\triangle_0 &:=& \{ (x_0, \ldots, x_n) \in \triangle: x_1+x_n >x_0\} \\
\triangle_1 &:=& \{ (x_0, \ldots, x_n) \in \triangle: x_1+x_n <x_0\}
\end{eqnarray*}
and define the slow-Triangle map
$T:\triangle_0\cup \triangle_1 \rightarrow \triangle $ by
\begin{eqnarray*}
T(x_0, \ldots, x_n) &=&  \left\{ \begin{array}{cc} T_0(x_0, \ldots, x_n),&\, \text{if } x_1+x_n>x_0 \\
                                                       T_1(x_0, \ldots, x_n), &\, \text{if } x_1+x_n<x_0 \end{array} \right.  \\
                                                       &=&  \left\{ \begin{array}{ccc} (x_1, x_2, \ldots , x_n, x_0-x_1), &\, \text{if }  x_1+x_n>x_0  \\
                                                       (x_0-x_n, x_1, x_2, \ldots, x_n), &\, \text{if }  x_1+x_n<x_0  \end{array} \right.  \\
\end{eqnarray*}
By writing the row vector $(x_0, \ldots, x_n)$ instead as a column vector, the action of $T$ is given by left multiplication by $(n+1)\times (n+1)$ matrices:

\begin{eqnarray*}
T \left( \begin{array}{c} x_0 \\ x_1 \\ \vdots  \\ x_n  \end{array}  \right)&=&  \left\{ \begin{array}{ccc} T_0 \left( \begin{array}{c} x_0 \\ x_1 \\ \vdots  \\ x_n  \end{array}  \right),&\, \text{if }  x_1+x_n>x_0 \\
                                                       T_1  \left( \begin{array}{c} x_0 \\ x_1 \\ \vdots  \\ x_n  \end{array}  \right),&\, \text{if } x_1+x_n<x_0 \end{array} \right.  \\
                                                       &=&  \left\{ \begin{array}{ccc}   \left( \begin{array}{c} x_1 \\ x_2 \\ \vdots \\ x_n  \\ x_0-x_1  \end{array}  \right),&\, \text{if } x_1+x_n>x_0 \\
                                                       \left( \begin{array}{c} x_0-x_n \\ x_1 \\ x_2 \\ \vdots  \\ x_n  \end{array}  \right), &\, \text{if } x_1+x_n<x_0\end{array} \right.  \\
  \end{eqnarray*}
where
\[
    T_0 = \begin{pmatrix} 0 & 1 & 0  & \cdots & 0\\ 0 & 0 & 1 & \cdots & 0 \\  && \vdots && \\0 & 0 & 0 & \cdots & 1  \\1 & -1 &
        0 & \cdots & 0 \end{pmatrix}
    \quad\text{and}\quad T_1= \begin{pmatrix} 1 & 0 & 0  & \cdots &  0 & -1\\ 0 & 1 &0&   \cdots &  0 & 0 \\ && \vdots &&&\\ 0 & 0 & 0 & \cdots & 0 & 1  \\\end{pmatrix}\]
Thus for $n=2$, we have 
\[
    T_0 = \begin{pmatrix} 0 & 1 & 0 \\ 0 & 0 & 1 \\ 1 & -1 &
        0 \end{pmatrix}
    \quad\text{and}\quad T_1= \begin{pmatrix} 1 & 0 & -1 \\ 0 & 1 & 0
        \\ 0&0&1 \end{pmatrix}
\]
For later use, note that 
\[
    \tau_0 \coloneqq T_0^{-1} =
    \begin{pmatrix}
        1 & 0 & 0&  \cdots & 0 &1 \\
        1 & 0 & 0 & \cdots & 0 & 0 \\ 0 &1&  0 &\cdots & 0 & 0 \\ 0&0&1& \cdots & 0 & 0\\
        & & \vdots &  &  & \\
        0 & 0 & 0 & \cdots & 1 & 0\end{pmatrix} \quad\text{and}\quad
    \tau_1\coloneqq T_1^{-1}= \begin{pmatrix} 1 & 0 & 0  & \cdots &  0 & 1\\ 0 & 1 &0&   \cdots &  0 & 0 \\ && \vdots &&&\\ 0 & 0 & 0 & \cdots & 0 & 1  \\\end{pmatrix}
\]
which for $n=2$ is
\[
    \tau_0 \coloneqq T_0^{-1} =
    \begin{pmatrix}
        1 & 0 & 1 \\
        1 & 0 & 0 \\ 0 & 1 & 0 \end{pmatrix} \quad\text{and}\quad
    \tau_1\coloneqq T_1^{-1}= \begin{pmatrix} 1 & 0 & 1 \\ 0 & 1 & 0
        \\ 0&0& 1 \end{pmatrix}.
\]

As a technical aside, it is important that both the matrices $\tau_0$
and $\tau_1$ have nonnegative entries, as it is in  the case for $\Phi_0$
and $\Phi_1$ in Section \ref{sec:Fmatrix}.  As we will see, this is
what allows the Triangle map to be used to understand partitions.
This is not necessarily the case for all multidimensional continued
fraction algorithms.

Following the construction in Section \ref{2 case}, we jump to the definition of the extended map given in Section \ref{sec:extended-farey}. The \emph{extended slow-Triangle map} $\tilde{T}$, which can be thought of as the natural extension of $T$, is defined by
\begin{eqnarray*}
\tilde{T} ( (n_1, \ldots , n_m) \times [k_1,\ldots , k_m]) &=&   \left\{ \begin{array}{cc} \tilde{T}_0 ( (n_1, \ldots , n_m) \times [k_1,\ldots , k_m]), &\, \text{if } n_2+n_m >n_1 \\
                                                        \tilde{T}_1 ( (n_1, \ldots , n_m) \times [k_1,\ldots , k_m]), &\, \text{if } n_2+n_m <n_1 \end{array} \right.  \\
                                                                                                              &=&  \left\{ \begin{array}{cc}    (n_2, n_3, \ldots, n_m, n_1-n_2), \times [k_1+ k_2, k_3, k_4, \ldots, k_m, k_1 ], &\, \\
                                                         \text{if } n_2+n_m >n_1 \\
                                                      (n_1-n_m, n_2, n_3, \ldots, n_m) \times [ k_1, \ldots, k_{m-1}, k_1 +  k_m  ], &\, \\
                                                       \text{if } n_2+n_m <n_1 \end{array} \right.  \\
\end{eqnarray*}
which can be read as the action of two $m\times m$ matrices on column
vectors in $\R^{2m}$, with the matrices
$$ \left( \begin{array}{cc} T_0 & 0  \\  0 & \tau_0^\top \end{array} \right) , \;   \left( \begin{array}{cc} T_1 & 0  \\  0 & \tau_1^\top \end{array} \right)\, .$$

\subsection{Link with integer partitions}

We now simply repeat what we did in Section \ref{2 case}, but now for
the slow-Triangle map. Consider a partition $(n_1^{k_1}, \ldots , n_m^{k_m}) \vdash n$ with $n_1>\cdots > n_m$, also
written as $(n_1, \ldots , n_m) \times [k_1,\ldots, k_m]\vdash n$. As
$n_1>\cdots  >n_m>0$, we can act on by the extended map. The key, both in
the definition of the natural extension and for our use in partition
theory, is that if
$$(n_1^{k_1}, \ldots , n_m^{k_m}) \vdash n,$$
then
$$ (n_2^{k_1+k_2} , n_3^{k_3}, \ldots ,n_m^{k_m}, (n_1-n_2)^{k_1}) \vdash n$$
and
$$( (n_1-n_m)^{k_1}, n_2^{k_2},\ldots, n_{m-1}^{k_m-1}, n_m^{k_1+k_m}) \vdash n.$$
This proves the following result.

\begin{prop} The extended slow-Triangle map $\tilde{T}$ sends a partition of $n$ to a new partition of $n$. Thus if $(n_1, \ldots , n_m) \times [k_1, \ldots , k_m]\vdash n,$ then
$$\tilde{T}((n_1, \ldots , n_m) \times [k_1, \ldots , k_m]) \vdash n.$$
\end{prop}

We can iterate the extended map $\tilde T$ and create an orbit of partitions.  As before, start with some
$$(n_1, \ldots , n_m)\times [1,0,\ldots, 0] \vdash n.$$
This a $0$th generation partition, which we write as
\[
    (n_1(0), \ldots , n_m(0)) \times [k_1(0),\ldots, k_m(0)].
\]
Acting on this vector by $\tilde{T}$ gets us the first generation
partition
\[
    \tilde{T}( (n_1(0), \ldots , n_m(0)) \times [k_1(0),\ldots , k_m(0)])
    = (n_1(1), \ldots , n_m(1)) \times [k_1(1),\ldots , k_m(1)].
\]
and recursively we obtain
$(n_1(a), \ldots , n_m(a)) \times [k_1(a),\ldots , k_m(a)]$ for  integers $a\ge 1$. 
  This relates the dynamics of the map $T$ with the sequence of partitions
obtained by $\tilde T$.

An example is :\label{orbit}
\[
    \begin{array}{c|c|c|c|c||c|c|c|c}
        a &( x, y )& n_1(a) & n_2(a) & n_3 (a)& k_1 (a)& k_2(a) & k_3(a)& \tilde{T}\\
        \hline
        0& (9/11, 4/11)&  11& 9 &  4    &  1 & 0  & 0 & \tilde{T}_0  \\
        1  &(4/9,  2/9)&     9    &   4    &  2  &  1&  0 &   1 &  \tilde{T}_1  \\
        2 &(4/7,  2/7)&       7  &   4    &   2 &  1& 0  &   2 & \tilde{T}_1  \\
        3  &(4/5, 2/5)&       5  &    4   &  2  &1  & 0  &   3 & \tilde{T}_0  \\
        4  &(2/4, 1/4)&       4  &   2    & 1   & 1 &  3 &  1  & \tilde{T}_1  \\
        5  &(2/3, 1/3)&       3  &   2    &  1  &  1&  3 &   2   \\
    \end{array}
\]

We stop here, for now, as $(2/3, 1/3)$ is on the line $x+y=1$. We will
deal with these boundary type points in Section \ref{$T_D$}.

Another example, in one dimension higher, is

\[
    \begin{array}{c|c|c|c|c|c||c|c|c|c|c}
      a   &( x, y, z )& n_1(a) & n_2(a) & n_3 (a)& n_4(a) & k_1 (a)& k_2(a) & k_3(a)& k_4(a) & \tilde{T}\\
        \hline
        0& (7/14, 6/14, 5/14)&  14 & 7 &  6    &  5 & 1&0  & 0&0 & \tilde{T}_1  \\
        1  &(6/9,  5/9, 2/9)&   9 & 7 &  6  &  5 & 1&0  & 0&1 & \tilde{T}_0 \\
        2 &(6/7,  5/7,2/7)&       7  &   6    & 5&  2 &  1& 0  &   1&1 & \tilde{T}_0  \\
        3  &(5/6, 2/6,1/6)&       6  &    5  &  2  &1  &1 &1& 1  &   1 & \\
        
    \end{array}
\]
Here we stop, as  $5/6 + 1/6= 1$ and is hence on the line $x+z=1$.  Again, in Section \ref{$T_D$} we will see how to continue this orbit.

In analogue to Proposition \ref{uniqueness} we have:

\begin{prop}
    For $n_2+n_m>n_1$, the partitions 
    \[ (n_1, \ldots , n_m) \times [1,0, \ldots , 0]\quad \text{and} \quad (n_1, n_3, \ldots, n_m, n_1-n_2) \times [1,0, \ldots , 0]\]
    will have the same generations after generation $0$.
\end{prop}

\begin{proof}
    As $n_2+n_m>n_1$, we first apply $\tilde{T}_0$ to
   $(n_1, \ldots , n_m) \times [1,0, \ldots , 0]$ and $\tilde{T}_1$ to \linebreak
   $(n_1, n_3, \ldots, n_m, n_1-n_2) \times  [1,0, \ldots , 0]$ to get 
    \begin{eqnarray*}
        \tilde{T}_0((n_1, \ldots , n_m) \times [1,0, \ldots , 0] &=&  (n_2, n_3,\ldots, n_m, n_1-n_2) \times [1,0, \ldots, 0, 1] \\
        &=& \tilde{T}_1((n_1, n_3, \ldots, n_m, n_1-n_2) \times  [1,0, \ldots , 0]).
    \end{eqnarray*}
    Thus after the first application of the extended slow-Triangle
    map, we have the same partition.
\end{proof}

\subsection{Conjugation and Palindromes for the slow-additive-Triangle map}

Here we extend the results of Section \ref{palindrome 2 dim} to the
slow-Triangle map. A partition
$(n_1, \ldots , n_m) \times [k_1\ldots , k_m]$ is conjugate to
$$(k_1+\ldots + k_m, k_1+\ldots + k_{m-1} , \ldots k_1 ) \times [n_m, n_{m-1}-n_m, \ldots, n_1-n_2]$$ (see Section \ref{partition}). We again have that the extended slow-Triangle map reflects and respects conjugation.

\begin{prop}\label{conjugation} The diagram
    \begin{eqnarray*}
       (n_1, \ldots , n_m) \times [k_1\ldots , k_m]  & \sim_{\mathcal{C}} &(k_1+\ldots + k_m, \ldots k_1 ) \times [n_m, n_{m-1}-n_m, \ldots, n_1-n_2]\\
        \tilde{T}_0 \downmapsto &   &   \upmapsto\tilde{T}_0\\
        (n_2,\ldots, n_1-n_2) \times [k_1+k_2, k_3, \ldots, k_m, k_1] & \sim_{\mathcal{C}} &
        (2k_1+k_2+ \ldots+ k_m, k_1+k_2+ \ldots +k_m, k_1+k_2) \\
        &&\times [n_1-n_2, n_{m-1}+ n_m-n_1,\ldots , n_2-n_3]
\end{eqnarray*}
when $n_2+ n_m > n_1 $ and the diagram
\begin{eqnarray*} (n_1, \ldots , n_m)  \times [k_1\ldots , k_m]  & \sim_{\mathcal{C}} &(k_1+\ldots + k_m, \ldots k_1 ) \times [n_m, n_{m-1}-n_m, \ldots, n_1-n_2] \\
\tilde{T}_1 \downmapsto &   &   \upmapsto \tilde{T}_1\\
(n_1-n_m, n_2, \ldots , n_m) \times [k_1,\ldots, k_{m-1}, k_1+ k_m  ] &  \sim_{\mathcal{C}} & (2k_1+k_2+ \ldots +k_m, k_1+\ldots + k_{m-1},\ldots ,  k_1)\\
&&  \times [n_m, n_{m-1}-n_m, \ldots, n_1-n_m- n_2]
\end{eqnarray*}
when $ n_2+ n_m <  n_1  $ are both commutative.
\end{prop}

As before, the proof is a simple calculation. Thus for the
slow-Triangle map we still have that the act conjugation can be viewed
as simply reversing the arrows of the extended map, giving us the
following result.

\begin{theorem}[The Power of Palindromes: Higher Dimension]
    Suppose that $ \tilde{T}_{i_1} , \ldots , \tilde{T}_{i_N}$ is a
    sequence of extended slow-Triangle maps such that
    $$ \tilde{T}_{i_N}\circ \cdots \circ  \tilde{T}_{i_1}(  (n_1, \ldots , n_m) \times   [k_1\ldots , k_m] )=   (\bar{n}_1, \ldots , \bar{n}_m) \times   [\bar{k}_1\ldots , \bar{k}_m]    $$
    Then
    $$\begin{array}{c} \tilde{T}_{i_1}\circ \cdots \circ  \tilde{T}_{i_N}((\bar{k}_1+\ldots + \bar{k}_m, \ldots \bar{k}_1 ) \times [\bar{n}_m, \bar{n}_{m-1}-\bar{n}_m, \ldots, \bar{n}_1-\bar{n}_2] )\\
        = (k_1+\ldots + k_m, \ldots k_1 ) \times [n_m,  n_{m-1}- n_m, \ldots, n_1- n_2]\end{array}$$

\end{theorem}


\subsection{On generalizing Theorem \ref{2parts-farey} and Theorem \ref{two sum}}

A difficulty arises when one wants to show the analogue of
Theorem \ref{2parts-farey} and of Theorem \ref{two sum}. When first
starting to explore the link between the slow-Triangle map and
partitions, we of course did many examples.  We set up charts of all
possible orbits for various $(n_1, n_2, n_3) \times [1,0,0]$ and we
quickly realized that the multiplicity vectors $[k_1,k_2, k_3]$ in
these orbits could only fit certain patterns.  These difficulties remain in higher dimensions.

\begin{definition} A vector of multiplicities $[k_1, \ldots , k_m] $ is
    \emph{allowable} if there are integers $n_1>\cdots >n_m >0$ with
    $n_1\not= n_2+n_m$ and integers $a_1>\cdots >a_m >0$ such that the
    partition $(a_1,\ldots , a_m) \times [k_1, \ldots, k_m]$ is a descendant
    of $(n_1,\ldots , n_m) \times [1,0, \ldots ,0]$.
\label{def-allowable}
\end{definition}
This leads to the natural question of which partitions into $m$ parts stem from an iteration of the extended slow-Triangle map starting from a root $(n_1,\ldots , n_m) \times [1,0, \ldots ,0]$.
What matters here, as we will see, is the multiplicity vector $[k_1, \ldots, k_m]$.
Now under $\tilde{T}_0$ 
we have 
$$[k_1, \ldots, k_m] \rightarrow [k_1+k_2, k_3, \ldots, k_m, k_1]$$
and 
under $\tilde{T}_1$ 
we have 
$$[k_1, \ldots, k_m] \rightarrow [k_1,k_2, \ldots, k_{m-1}, k_m+  k_1].$$

We need a few technical lemmas.  For notation, suppose we are iterating a partition.  At the $pth$ step, label the corresponding multiplicity as 
$$[k_1(p), \ldots, k_m(p)].$$  Further, recall that 
$$[  k_1(0), \ldots, k_m(0)] = [1,0, \ldots ,0].$$

By simply looking at the maps $\tilde{T}_0$  and $\tilde{T}_1$, we can see 
\begin{lemma}
For all $p$ we have $k_1(p) >0.$
\end{lemma}

Noting that $\tilde{T}_0$ sends the first term in the multiplicity to the last and shifts to the left the others, and that  $\tilde{T}_1$  leave all the terms alone, save for adding the first term to the last, we see that 

\begin{lemma}
If there is an $i>1$ so that 
$$k_i(p), k_{i+1}(p), \ldots, k_m(p) >0$$
then 
$$k_i(p+1), k_{i+1}(p+2), \ldots, k_m(p+1) >0$$

\end{lemma}
  Further
\begin{lemma}
If $k_1(p) < k_1(p+1)$, then we must have $k_2(p) >0.$

\end{lemma}

Putting these lemmas together gives us 
\begin{lemma}
A vector $ [k_1, 0, k_3, \ldots, k_m] $ with $k_1>1$ is not an  allowable multiplicity.
\end{lemma}

\begin{proof}
This follows since the initial first term must be $1$.
\end{proof}

\begin{lemma}
A vector $ [k_1, \ldots, k_{m-1}, 0] $ with $k_2>0,$ is not an allowable multiplicity.
\end{lemma}

\begin{proof}
Applying either   $\tilde{T}_0$  or $\tilde{T}_1$ to the initial $[1,0, \ldots ,0]$ will give us 
$$[1,0, \ldots ,1]$$
The last term can never return to zero.  \end{proof}

\begin{lemma} \label{nonallowable}
A vector $ [k, k_2, \ldots, k_{m-1}, k ] $ with $k>0$ is not an  allowable multiplicity.
\end{lemma}
\begin{proof}
Suppose that $ [k, k_2, \ldots, k_{m-1}, k, ] $  is allowable.  Then it be the image of an allowable multiplicity from $\tilde{T}_0$  or from $\tilde{T}_1$.  Now, the inverse of this multiplicity by $\tilde{T}_0$ is 
$$[k, 0, k_2, \ldots , k_{m-1}]$$
which is not allowable.
The inverse of this multiplicity by $\tilde{T}_1$ is 
$$[k, k_2,  \ldots , k_{m-1}, 0]$$
which is  also not allowable.

\end{proof}

As an example, this means that the partition of 20 given by 
$$(5,4,3) \times [2,1,2]$$
cannot be obtained from any possible
sequences of $\tilde{T}_0$  and $\tilde{T}_1$ from any 
$$(20, n_2, n_3) \times [1,0,0].$$

The goal of the next section is to get around this difficulty.

\section{The map $\tilde{T}_D$ on the boundary} \label{$T_D$}
\subsection{The definition and naturalness of $\tilde{T}_D$ }

The original additive-slow-Triangle map $T$ acting on the space $x_0>x_1>\cdots >x_m>0$ is  simply not defined on the hyperplane $x_0=x_1+x_m$.  Before now, much of the work on this map has been concerned with it as a dynamical system.  As this hyperplane has measure zero, it could be conveniently ignored. But when applied to partitions, this means we would be ignoring many perfectly reasonable partitions.  In this section we will show how to extend the map $\tilde{T}$ to the hyperplane $x_0=x_1+x_m$. (This map was originally defined in work done in parallel to this paper in \cite{Baalbaki-Garrity}.)

We start with an example.  Consider the partition $(6,5,2,1) \times [1,1,1,1]\vdash 14.$  If we naively apply $\tilde{T}_0$ and $\tilde{T}_0$, we have

\begin{eqnarray*}
(6,5,2,1) \times [1,1,1,1] & \xrightarrow{\tilde{T}_0} & (5,2,1,1) \times [2,1,1,1]\\
(6,5,2,1) \times [1,1,1,1] & \xrightarrow{\tilde{T}_1} & (5,5,2,1) \times [1,1,1,2]\\
\end{eqnarray*}
Both of these images would be a strange way for writing the partitions, as the natural notation would be to concatenate common terms.  But after concatenation, we get the same partition: 
$$(6,5,2,1) \times [1,1,1,1] \rightarrow (5,2,1) \times [2,1,2].$$
This leads to the following map $\tilde{T}_D$:
$$\tilde{T}_D( (n_1, \ldots, n_m) \times[k_1, \ldots, k_m] )= (n_2, \ldots, n_m) \times[k_1+k_2, k_3 , \ldots, k_{m-1}, k_1+ k_m]$$
whenever
$$n_1=n_2+n_m.$$

This allows us to extend the orbits in  \ref{orbit}:

\[
    \begin{array}{c|c|c|c|c||c|c|c|c}
        a &( x, y )& n_1(a) & n_2(a) & n_3 (a)& k_1 (a)& k_2(a) & k_3(a)& \tilde{T}\\
        \hline
        0& (9/11, 4/11)&  11& 9 &  4    &  1 & 0  & 0 & \tilde{T}_0  \\
        1  &(4/9,  2/9)&     9    &   4    &  2  &  1&  0 &   1 &  \tilde{T}_1  \\
        2 &(4/7,  2/7)&       7  &   4    &   2 &  1& 0  &   2 & \tilde{T}_1  \\
        3  &(4/5, 2/5)&       5  &    4   &  2  &1  & 0  &   3 & \tilde{T}_0  \\
        4  &(2/4, 1/4)&       4  &   2    & 1   & 1 &  3 &  1  & \tilde{T}_1  \\
        5  &(2/3, 1/3)&       3  &   2    &  1  &  1&  3 &   2  & \tilde{T}_D\\
        6& (1/2) & 2 & 1 & & 4& 3& 
    \end{array}
\]

and

\[
    \begin{array}{c|c|c|c|c|c||c|c|c|c|c}
      a   &( x, y, z )& n_1(a) & n_2(a) & n_3 (a)& n_4(a) & k_1 (a)& k_2(a) & k_3(a)& k_4(a) & \tilde{T}\\
        \hline
        0& (7/14, 6/14, 5/14)&  14& 7 &  6    &  5 & 1&0  & 0&0 & \tilde{T}_1  \\
        1  &(7/9,  6/9, 5/9)&   9& 7 &  6    &  5 & 1&0  & 0&1 & \tilde{T}_0 \\
        2 &(6/7,  5/7,2/7)&       7  &   6    & 5&  2 &  1& 0  &   1&1 & \tilde{T}_0  \\
        3  &(5/6, 2/6,1/6)&       6  &    5  &  2  &1  &1 &1& 1  &   1 &  \tilde{T}_D  \\
        4 & (2/5, 1/5) & 5 & 2 & 1 & & 2 & 1 & 2 & &  \tilde{T}_1  \\
      5 & (2/4, 1/4) & 4 & 2 & 1 & & 2 & 1 & 4 & &  \tilde{T}_1  \\
 6 & (2/3, 1/3) & 3 & 2 & 1 & & 2 & 1 & 6 & &  \tilde{T}_D  \\
  7 & (1/2) & 2 & 1 &  & & 2 & 5 &  & &   \\
    \end{array}
\]
Thus changing the dimension allows us to create longer trees.  

\subsection{Capturing non-allowable partitions}\label{Graphs}
We just saw that the partition 
$(5,2,1) \times [2,1,2]\vdash 14$ is the image of $(6,5,2,1)\times [1,1,1,1]$ under the map $\tilde{T}_D$.  But Lemma \ref{nonallowable} seems to state that $(5,2,1) \times [2,1,2]$ is nonallowable.  Of course this is not at all a contradiction, as Lemma \ref{nonallowable} is only about possible images of $\tilde{T}_0$ and $\tilde{T}_1$.  Thus allowing us to ``change dimensions'' allows us to capture all partitions.

Part of this comes down to the inverses of the three maps $\tilde{T}_0$, $\tilde{T}_1$, which are  $\tilde{T}_D$.  As described in \cite{Baalbaki-Garrity}, the maps $\tilde{T}_0$ and $\tilde{T}_1$ are one-to-one, each having the following inverses:
 for $k_1>k_m$ and $m>2$,
$$(n_1, \ldots , n_m) \times [k_1, \ldots, k_m] \xrightarrow{\tilde{T}_0^{-1}  }(n_1+ n_m, n_1, n_2, n_3,    \ldots, n_{m-1}) \times [k_m,k_1-k_m, k_2,  k_3, \ldots, k_{m-1}] $$
for $k_1<k_m$ and $m>2$,

$$(n_1, \ldots , n_m) \times [k_1, \ldots, k_m] \xrightarrow{\tilde{T}_1^{-1}  }   (n_1+ n_m, n_2,  \ldots ,    \ldots, n_{m}) \times [k_1, \ldots, \ldots , k_{m-1}, k_m-k_1],       $$
 for $k_1>k_2$ and $m=2,$
   $$ (n_1, n_2) \times [k_1, k_2] 
\xrightarrow{ T_0^{-1}  }   
   (n_1+ n_2, n_1) \times [k_2,  k_1-k_2]     $$
   and for 
$k_1<k_2$ and $m=2$
   $$ (n_1, n_2) \times [k_1, k_2] \xrightarrow{ T_1^{-1}  }  
   (n_1+ n_2, n_2) \times [k_1,  k_2-k_1]  .$$
   The map $\tilde{T}_D$ that changes dimensions though is not one-to-one. 
If $$K= \min  (k_1+k_2 , k_1+k_m ) - 1, $$
  then $T_D$ is a $K$ to one map. Explicitly, for any $k=1, \ldots, \min(k_1, k_{m-1}) -1,$ an inverse is 
 $$(n_1, \ldots , n_m) \times [k_1, \ldots, k_m] \xrightarrow{\tilde{T}_D^{-1}  } (n_1+n_m, n_1, \ldots, n_{m-1} ) \times [k, k_1-k, k_2, \ldots, k_{m-2}, k_{m-1} -k].$$
  
  This allows us to fit any partition into an eventual image of some 
  $(n, n_2, \ldots, n_m)\times [1, \ldots,0 ].$
  This creates a network of interrelations among partitions, a network that overall remains to be fully understood.

\section{On other multidimensional continued fraction algorithms}\label{otherMCFs}

There are many different multidimensional continued fraction algorithms. (To get a feel of how many there are, see Schweiger \cite{Schweiger4}  and Karpenkov \cite{Karpenkov13}.)  Each has their own strengths and weaknesses.  All are trying to generalize the many wonderful properties of traditional continued fractions to higher dimensional analogs.  Historically, the two main sources of inspiration have been trying to find good Diophantine approximation properties and trying to find methods for understanding algebraic numbers via periodicity properties (generalizing the classical fact that a real number is a quadratic irrational if and only if its continued fraction expansion is eventually periodic.) 

When we started this project, we assumed that each of the other well known multidimensional continued fraction algorithms would provide their own dynamical interpretation of partitions.  This is quite false, as we will now see.  To give a flavor of other multidimensional continued fractions algorithms and why they are not useful at all for studying partition numbers, we will look at the M\"{o}nkemeyer algorithm and then the Cassaigne algorithm.  We then turn to the language of triangle partition maps, which is an attempt to put various multidimensional continued fraction algorithms into a single framework.  We will see in terms of this framework, hardly any multidimensional continued fraction algorithm can be used in partition theory.  Thus in the context of linking dynamics with partition numbers, it seems that the triangle map is unusual (though as we will also see not unique).

\subsection{M\"{o}nkemeyer}
The M\"{o}nkemeyer map is a particularly good multidimensional continued fraction algorithm for generalizing the classical Minkowksi Question Mark function, as shown by Panti \cite{Panti}.  For ease of notation, we will only treat the case of $m=3.$

Here we will simply write down the map.
Set
\begin{eqnarray*}
\triangle &=& \{ (x,y) \in \R^2: 1 >x>y>0\} \\
\triangle_0 &=& \{ (x,y) \in \triangle: x+y >1\} \\
\triangle_1 &=& \{ (x,y) \in \triangle: x+y <1\}
\end{eqnarray*}
exactly that same as in Section  \ref{npartitions}.
The slow-M\"{o}nkemeyer map
$M:\triangle_0 \cup \triangle_1 \rightarrow \triangle $ is
\begin{eqnarray*}
M(x,y) &=&  \left\{ \begin{array}{cc} M_0(x,y), &\, \text{if } x+y>1 \\
                                                       M_1(x,y), &\, \text{if } x+y<1 \end{array} \right.  \\
                                                       &=&  \left\{ \begin{array}{ccc} \left( \frac{1-y}{x}, \frac{x-y}{x}  \right), &\, \text{if } x+y>1 \\
                                                        \left( \frac{x}{1-y}, \frac{x-y}{1-y}  \right), &\, \text{if } x+y<1 \end{array} \right.  \\
\end{eqnarray*}

As in  Section \ref{sec:Fmatrix} and Section \ref{npartitions},  we  pass from points $(x,y) $ in $\R^2$ to vectors $(z,x,y)$
in $\R^3$ (or points $(z:x:y)$ in $\R\P^2$) via sending $(x,y) $ to
$(1,x,y)$ with inverse map $(z,x,y)\rightarrow (x/z, y/z).$ Then, by
an abuse of notation, as before, we again set 
\begin{eqnarray*}
\triangle &:=& \{ (z,x,y) \in \R^3: z >x>y>0\} \\
\triangle_0 &:=& \{ (z,x,y) \in \triangle: x+y >z\} \\
\triangle_1 &:=& \{ (z,x,y) \in \triangle: x+y <z\}.
\end{eqnarray*}
This allows us to define the slow-M\"{o}nkemeyer map
$T:\triangle_0\cup \triangle_1 \rightarrow \triangle $ by
\begin{eqnarray*}
M(z, x,y) &=&  \left\{ \begin{array}{cc} M_0(z, x,y), &\, \text{if } x+y>z \\
                                                       M_1(z, x,y), &\, \text{if } x+y<z \end{array} \right.  \\
                                                       &=&  \left\{ \begin{array}{ccc} (x,z- y, x-y), &\, \text{if } x+y>z \\
                                                       (z-y, x, x-y), &\, \text{if } x+y<z \end{array} \right.  \\
\end{eqnarray*}
By writing the row vector $(z,x,y)$ as a column vector, the action of $M$ is given by left multiplication by $3\times 3$ matrices:

\begin{eqnarray*}
M \left( \begin{array}{c} z \\ x  \\ y  \end{array}  \right)&=&  \left\{ \begin{array}{ccc} M_0  \left( \begin{array}{c} z \\ x  \\ y  \end{array}  \right), &\, \text{if }  x+y>z \\
                                                       M_1  \left( \begin{array}{c} z \\ x  \\ y  \end{array}  \right), &\, \text{if } x+y<z \end{array} \right.  \\
                                                       &=&  \left\{ \begin{array}{ccc}   \left( \begin{array}{c} x \\ z-y  \\ x-y  \end{array}  \right), &\, \text{if } x+y>z \\
                                                       \left( \begin{array}{c} z-y \\ x  \\ x-y  \end{array}  \right), &\, \text{if } x+y<z \end{array} \right.  \\
  \end{eqnarray*}
where
\[
    M_0 = \begin{pmatrix} 0 & 1 & 0 \\ 1 & 0 & -1 \\ 0 & 1 &
        -1 \end{pmatrix}
    \quad\text{and}\quad M_1= \begin{pmatrix} 1 & 0 & -1 \\ 0 & 1 & 0
        \\ 0&1&-1 \end{pmatrix}
\]
and
\[
    m_0 \coloneqq M_0^{-1} =
    \begin{pmatrix}
        1 & 1 &- 1 \\
        1 & 0 & 0 \\ 1 & 0 & -1 \end{pmatrix} \quad\text{and}\quad
   m_1\coloneqq M_1^{-1}= \begin{pmatrix} 1 & 1 &- 1 \\ 0 & 1 & 0
        \\ 0&1&- 1 \end{pmatrix}.
\]
Note that the entries for both $m_0$ and $m_1$ have negative entries.  This is the technical reason why the M\"{o}nkemeyer map will not be good to understand partitions.

As seen in the earlier sections, we need to look at the \emph{extended slow-M\"{o}nkemeyer  map} $\tilde{M}$ (the natural extension of $M$),  which is
\begin{eqnarray*}
\tilde{M} ( (n_1, n_2, n_3) \times [k_1, k_2, k_3]) &=&   \left\{ \begin{array}{cc} \tilde{M}_0 ( (n_1, n_2, n_3) \times [k_1, k_2, k_3]), &\, \text{if } n_2+n_3 >n_1 \\
                                                        \tilde{M}_1 ( (n_1, n_2, n_3) \times [k_1, k_2, k_3]), &\, \text{if } n_2+n_3 <n_1 \end{array} \right.  \\
                                                                                                              &=&  \left\{ \begin{array}{cc}    (n_2, n_1-n_3, n_2-n_3) \times [    k_1+ k_2+k_3, k_1, -k_1 -k_3], &\, \text{if } n_2+n_3 >n_1\\
                                                      (n_1-n_3, n_2, n_2-n_3) \times [ k_1, k_1+k_2+ k_3, -k_1 -  k_3  ], &\, \text{if } n_2+n_3 <n_1 \end{array} \right.  \\
\end{eqnarray*}
which can be read as the action of two $6\times 6$ matrices on column
vectors in $\R^6$, with the matrices
$$ \left( \begin{array}{cc} M_0 & 0  \\  0 & m_0^\top \end{array} \right) , \;   \left( \begin{array}{cc} M_1 & 0  \\  0 & m_1^\top \end{array} \right)\, .$$

If we try to see, for example, what the partition 
$$(7, 5, 4) \times [3, 2, 4]\vdash 47$$
would map to under $\tilde{M}$, we get 
$$\tilde{M}((7, 5, 4) \times [3, 2, 4])=(5,3,2)\times [9, 3, -7].$$
That $-7$ for one of the multiplicities means that this dynamical system will not generate partitions.

\subsection{Cassaigne}
This algorithm is  of  fairly recent vintage.  Back in the early 1940s, Morse and Hedlund \cite{Marse-Hedlund} started investigating infinite words.  They showed that any infinite word made up from an alphabet of two letters whose linear complexity is bounded by $n$ (meaning that there are no more than $n$  subwords of length $n$) must actually be eventually periodic. An infinite word $w$  whose linear complexity is exactly $n+1$ is called Sturmian, meaning that there are exactly $n+1$  subwords of length $n$ in $w$.  Morse and Hedlund not only showed that Sturmian words exist but more so that all such words stem from continued fraction expansions of real numbers.  This is quite amazing.  (For more see Arnoux's work in Chapter 6 of N. Pytheas Fogg 's {\it Substitutions in Dynamics, Arithmetics, and Combinatorics} \cite{Fogg}.

The natural question then arises for methods to generalize this link of continued fractions with infinite words of low complexity. 
The Cassaigne algorithm, as described by Cassaigne,  Labb\'{e} and  Leroy \cite{Cassaigne}, is a multidimensional continued fraction algorithm that produces infinite words on three letters whose linear complexity is exactly $2n+1$.  
We could now describe how this algorithm acts on vectors in the cone $\triangle$ via matrix multiplication, etc.  We will instead simply write down 
 the \emph{extended slow-Cassaigne  map} $\tilde{M}$ ,  which is 
\begin{eqnarray*}
\tilde{C} ( (n_1, n_2, n_3) \times [k_1, k_2, k_3])              &=&  \left\{ \begin{array}{cc}    (n_2, n_3, n_2+n_3- n_1) \times [    k_1+ k_2, k_1+k_3, -k_1 ], &\, \text{if } n_2+n_3 >n_1\\
                                                      (n_1-n_3, n_2, n_2-n_3) \times [ k_1, k_1+k_2+ k_3, -k_1 -  k_3  ], &\, \text{if } n_2+n_3 <n_1 \end{array} \right.  \\
\end{eqnarray*}

Looking at the partition
$$(7, 5, 4) \times [3, 2, 4]\vdash 47$$
again,  we see that 
$$\tilde{C}((7, 5, 4) \times [3, 2, 4])=(5,4,2)\times [5, 7, -3].$$
That $-3$ for one of the multiplicities means that this dynamical system will  also not generate partitions.

\subsection{Triangle Partition Maps}

There are many multidimensional continued fractions.  For each of these we can look at the corresponding extended map and see if the maps on the multiplicities keeps all the terms positive.  This seems like an almost endless process.
In \cite{tripmaps}, a family of multidimensional continued fractions was created which seems to capture all known multidimensional continued fractions (in a certain well-defined sense).  The method is  to first create a list of $216$
different multidimensional continued fractions, parameterized by $S_3\times S_3 \times S_3$, where $S_3$ is the set of permutations of three elements. Thus given
some $(\sigma, \tau_0, \tau_1) \in S_3\times S_3 \times S_3$ we associate a multidimensional continued fraction which we will denote by $T(\sigma, \tau_0, \tau_1) $, which in turn stems from two maps
$$T_0(\sigma, \tau_0, \tau_1),\, T_1(\sigma, \tau_0, \tau_1) $$
 and corresponding extended versions
$$ \tilde{T}(\sigma, \tau_0, \tau_1),\, \tilde{T}_0(\sigma, \tau_0, \tau_1),\, \tilde{T}_1(\sigma, \tau_0, \tau_1).$$
It can be calculated that the Triangle map is $T(e,e,e)$, that M\"{o}nkemeyer is $T(e,132,23)$ and that Cassaigne is $T(e,23,23).$  

 For the generating list of 216 different multidimensional continued fraction algorithms, we can explicitly compute for each the corresponding extended version.  What struck us initially as  a somewhat surprising fact is that only four of these maps can be used to study partitions, meaning that for all but four some of the multiplicities will be negative.  These four are
 $$T(e,e,e),\, T(13,12,12),\, T(12,e,12),\, T(132, 12,e).$$
 Though the first is the Triangle map, the other three do not have names and have never really been studied  before.
 
 While all four of these maps give rise to different dynamical systems, there is the phenomenon of ``twinning'' (see  Section 8.1 in \cite{Garrity-McDonald}), which is that 
 $$\begin{array}{cc}
 
 \begin{array}{ccc}   
 \tilde{T}_0(e,e,e)  &=& \tilde{T}_1(13,12,12),\\
 \tilde{T}_1(e,e,e)  &=& \tilde{T}_0(13,12,12),
 \end{array} &   
 \begin{array}{ccc}   
 \tilde{T}_0(12,e,12)  &=& \tilde{T}_1(132, 12,e)\\
 \tilde{T}_1(12,e,12)  &=& \tilde{T}_0((132, 12,e)
 \end{array} 
 \end{array} $$
This means really there are only two different $T(\sigma, \tau_0,\tau_1)$ maps that can be used to study partitions.
As $ T(12,e,12)$ has never really been studied before, we will  briefly describe it here. Set
\begin{eqnarray*}
\triangle &:=& \{ (x,y) \in \R^2: 1 >x>y>0\} \\
\triangle_0(12,e,12) &:=& \{ (x,y) \in \triangle:2y >x\} \\
\triangle_1(12,e,12) &:=& \{ (x,y) \in \triangle: 2y <x\}
\end{eqnarray*}

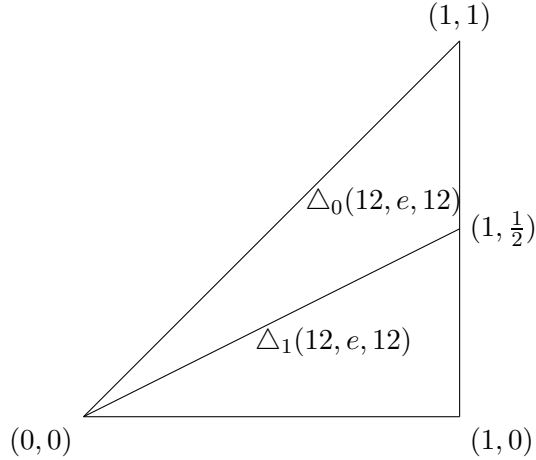
\begin{center}
\begin{tikzpicture}[scale=5]
\draw(0,0)--(1,0);
\draw(0,0)--(1,1);
\draw(1,0)--(1,1);
\draw(0,0)--(1,1/2);
\node[] at (4/5,4/7){$\triangle_0(12,e,12) $};
\node[] at (4/6,1/5){$\triangle_1(12,e,12) $};
\node[below left]at(0,0){$(0,0)$};
\node[below right]at(1,0){$(1,0)$};
\node[above]at(1,1){$(1,1)$};
\node[right] at (1,1/2){$(1,\frac 12)$};
\end{tikzpicture}
\end{center}

The  map
$T(12,e,12):\triangle_0(12,e,12) \cup \triangle_1(12,e,12) \rightarrow \triangle $ is
\begin{eqnarray*}
T(x,y) &=&  \left\{ \begin{array}{cc} T_0(x,y), &\, \text{if } 2y>x \\
                                                       T_1(x,y), &\, \text{if } 2y<x \end{array} \right.  \\
                                                       &=&  \left\{ \begin{array}{ccc} \left( \frac{y}{1+y-x}, \frac{x-y}{1+y-x}  \right), &\, \text{if } 2y>x \\
                                                        \left( \frac{x-y}{1-y}, \frac{y}{1-y}  \right), &\, \text{if } 2y<x \end{array} \right.  \\
\end{eqnarray*}

The  extended  $\tilde{T}(12,e,12) $, which can be thought of as the natural extension of $T(12,e,12)$, is 
\begin{eqnarray*}
\tilde{T}(12,e,12) ( (n_1, n_2, n_3) \times [k_1, k_2, k_3]) &=&   \left\{ \begin{array}{cc} \tilde{T}_0(12,e,12) ( (n_1, n_2, n_3) \times [k_1, k_2, k_3]), &\, \text{if } 2n_3 >n_2 \\
                                                        \tilde{T}_1(12,e,12) ( (n_1, n_2, n_3) \times [k_1, k_2, k_3]), &\, \text{if } 2n_3 <n_2 \end{array} \right.  \\
                                                                                                              &=&  \left\{ \begin{array}{cc}    (n_1+n_3-n_2, n_3, n_2-n_3) \times [    k_1, k_2+ k_3, k_1+k_2 ], &\, \text{if } 2n_3 >n_2\\
                                                      (n_1-n_3, n_2-n_3, n_3) \times [ k_1, k_2, k_1 + k_2+ k_3  ], &\, \text{if } 2n_3 <n_2 \end{array} \right.  \\
\end{eqnarray*}
As the multiplicities will never be negative, we have orbits of partition numbers.    But these orbits will be different than the corresponding orbits for the Triangle map, as can be seen by comparing the orbit of 
$$(11,9,4) \times [1,0,0]$$
via $\tilde{T}(12,e,12) $, which is
\[
    \begin{array}{c|c|c|c|c||c|c|c|c}
        m &( x, y )& n_1(m) & n_2(m) & n_3 (m)& k_1 (m)& k_2(m) & k_3(m)&\tilde{T}_i(12,e,12) \\
        \hline
        0& (9/11, 4/11)&  11& 9 &  4    &  1 & 0  & 0  &\tilde{T}_1(12,e,12) \\
        1  &  (4/5, 1/5)  &7  &    5      &   4     &   1  &  0 &  1  &  \tilde{T}_0(12,e,12)     \\
        2 &  (4/6, 1/6) &6  &      4    &   1     &   1  & 1  &  1  &   \tilde{T}_1(12,e,12)    \\

        3  &(3/5, 1/5)  & 5  &     3     &    1    & 1   &  1 &  3  &    \tilde{T}_1(12,e,12)   \\

        4  & (2/4, 1/4) & 4  &    2      &     1   &    1 &1   & 5   &       \\

    \end{array}
\]
(stopping  here, for now, as $(2/4, 1/4)$ is on the line $2y=x$) 
to the orbit under the triangle map $\tilde{T}(e,e,e) $ given in the table  in Section \ref{orbit}.

Unfortunately, while $\tilde{T}(12,e,12) $ does provide a map of partitions to partitions, almost any example will show that this map does not respect conjugation, and hence the analog of Proposition \ref{conjugation} is false.

\subsection{Selmer and Brun}
Recently Matthew Phang \cite{Phang}  has shown that two classical multi-dimensional continued fractions, the Selmer algorithm (see Chapter 7 in Schweiger \cite{Schweiger4} ) and the Brun algorithm (Chapter 4 in \cite{Schweiger4} ), both can be used to produce maps from partitions to partitions.  Neither though respect conjugation and hence for both of these the the analog of Proposition \ref{conjugation}  are also false.

\section{Questions}

We view this paper as only a start.  For example, in \cite{Baalbaki-Garrity} many new identities among partitions are given using the dynamics of the triangle map applied to partitions.  This work also gives a process to discover many new partition identities.  This strikes as quite promising.

In this paper, we have spent a lot of time on the special case of partitioning a number into two numbers, with multiplicity, in large part due to the richness behind the classical Farey map and its corresponding Farey tree.   This correspondence is what is   key to our formula for $p(2,n)$.   Recently, two of the authors, with Sara Munday, \cite{Bonanno- Del Vigna-Munday} have developed an analogous tree structure for the slow-Triangle map (see also \cite{Bonanno-Del Vigna}).  Thus suggests that there will be a rich analog of the Farey tree, a possible Triangle graph.  This is underlying Subsection \ref{Graphs}.  We hope to pursue this in future work.

\section*{Acknowledgements}

We would like to thank George Andrews, Tim Huber, Lori Pedersen, Jimmy McLaughlin and  Ken Ono for help. Claudio Bonanno, Alessio Del Vigna and Stefano Isola are partially supported by the PRIN Grant 2017S35EHN of the Ministry of Education, University and Research (MIUR), Italy. The research of Claudio Bonanno, Alessio Del Vigna and Stefano Isola is also part of their activity within the UMI group ``DinAmicI'' (\texttt{www.dinamici.org}) and the Gruppo Nazionale di Fisica Matematica, INdAM, Italy.

\end{document}